\documentclass[reqno,10pt,a4paper]{amsart}
\usepackage{amssymb}
\usepackage{latexsym,amssymb,amsfonts}

\begin{document}
\author{Robert Sh. Liptser}
\author{Anatolii A. Pukhalskii}
\title[Limit theorems on large deviations for semimartingales]
{Limit theorems on large deviations for semimartingales}

\address{Institute for Information Transmission Sciences, 19,
Bolshoy Karetnii (Ermolovoy str.), 101447, Moscow, Russia}

\keywords{Large deviations, Semimartingale, Skorokhod space,
Cumulant, Fenchel-Legendre transform, Multiplicative
decomposition, Change of probability measure.}

\subjclass{60F10}

\maketitle
\begin{abstract}
We consider a sequence $X^n=(X^n_t)_{t\ge 0},n\ge 1$ of
semimartingales. Each $X^n$ is a weak solution to an It\^o
equation  with respect to a Wiener process and a Poissonian
martingale measure and is in general non-Markovian process. For
this sequence, we prove the large deviation principle in the
Skorokhod space $D=D_{[0,\infty)}$. We use a new approach based on
of exponential tightness. This allows  us to establish the large
deviation principle under weaker assumptions than before.
\end{abstract}

\renewcommand{\b}[1]{\mbox{\boldmath $#1$}}
\renewcommand{\frak}{\mathfrak}
\renewcommand{\Bbb}{\mathbb}
\newtheorem{problem}{Problem}
\newtheorem{lemma}{\quad\bf Lemma}[section]
\newtheorem{proposition}{Proposition}[section]
\newtheorem{corollary}{Corollary}
\newtheorem{remark}{\quad \it Remark}
\newtheorem{conjecture}{Conjecture}
\newtheorem{algorithm}{Algorithm}
\newtheorem{theorem}{\bf Theorem}[section]
\newtheorem{exercise}{Exercise}[section]
\numberwithin{equation}{section} \theoremstyle{plain}
\theoremstyle{definition}
\newtheorem{definition}{\sc Definition}[section]
\theoremstyle{remark} \newtheorem{example}{Example}[section]

\vskip 1in \noindent{Main notations}

$$
\begin{array}{ll}
(\Omega,\mathcal{F},\mathbb{F}=(\mathcal{F}_t)_{t\ge 0},P), &
\text{a stochastic basis};
 \\
T_L(\mathbb{F}), & \text{the set of stopping times (relative to a
filtration $\mathbb{F}$}
\\ & \text{not  exceeding $L$};
\\
D=D_{[0,\infty)}(D_{0,T]}), & \text{the Skorokhod space of all
right continuous,}
\\
& \text{having left hand limits real valued functions}
\\
& X=(X_t)_{t\ge 0}(X=(X_t)_{0\le t\le T});
\\
C=C_{[0,\infty)}(C_{0,T]}), & \text{the space of all right
continuous functions}
\\
&\text{from $D_{[0,\infty)}(D_{[0,T]})$;}
\\
(E,\mathcal{E})& \text{a Blackwell space};
\\
\mathcal{B}(R), \ \mathcal{B}(R_+), \ \mathcal{B}(R_0),& \text{the
Borel $\sigma$-fields on $R, \ R_+, \ R_0:=R\setminus\{0\}$};
\\
\mathbb{D}=(\mathcal{D})_{t\ge 0}, &\text{the family of
$\sigma$-algebras $\mathcal{D}_t=\sigma(X_s,s\le t), X\in D$;}
\\
\rho(.,.), & \text{the Lindvall-Skorokhod metric on $D$};
\\
``\xrightarrow{P}'', & \text{convergence in probability};
\\
x\wedge y=\min(x,y), & x\vee y=\max(x,y);
\end{array}
$$

$$
X^*_t=\sup_{s\le t}|X_s|, \ X^*_{t-}=\sup_{s<t}|X_s|, \ \triangle
X_t=X_t- X_{t-}, \ X\in D;
$$

\newpage
\section{\bf Setting, basic concept, method}
\label{sec-1}\noindent {\bf 1}. Let $X^n=(X^n_t)_{t\ge 0}, n\ge 1$
be a sequence of stochastic processes with path in $D$. For each
$n$, $X^n$ is a semimartingale on a stochastic basis
$(\Omega,\mathcal{F},\mathbb{F}^n=(\mathcal{F}^n_t)_{t\ge 0},P)$
and is a weak solution to the Ito equation \cite{Gsk}
\begin{multline}\label{1.1}
X^n_t=x+\int_0^ta(s,X^n)ds+\frac{1}{\sqrt{n}}\int_0^tb(s,X^n)dW^n_s
\\
+ \frac{1}{n}\int_0^t\int_Ef(s,X^n,u)[p^n-q^n](ds,du),
\end{multline}
where $x\in R$, $(W^n_t)_{t\ge 0}$ is a Winer process,
$p^n(dt,du)$ is an integer valued random measure on $(R_+\times
E,\mathcal{B}(R_+)\otimes\mathcal{E})$, $q^n(dt,du)$ is the
compensator of $p^n(dt,du)$ with
\begin{equation}\label{1.2}
q^n(dt,du)=ndtq(du),
\end{equation}
$q(du)$ is a measure on $(E,\mathcal{E})$ (for all the definitions
and facts from the martingale theory we refer the reader to
\cite{LSMar} and \cite{JS}.

We assume that functionals $a(t,X)$, $b(t,X)$ and $f(t,X,u)$
($t\in R_+,X \in D, u\in E$) are $\mathcal{P}(\mathbb{D})$- and
$\widetilde{\mathcal{P}} (\mathbb{D})$-measurable and satisfy the
following conditions:

{\bf I}. (linear growth). For all $X\in D$
\begin{eqnarray*}
|a(t,X)|&\le& l_t(1+X^*_{t-}),
\\
|b(t,X)|&\le& l_t(1+X^*_{t-})
\\
|f(t,X,u)|&\le& h(u)l_t(1+X^*_{t-})
\end{eqnarray*}
for almost all $t\in R_+$ with respect to the Lebesgue measure,
where $l_t$ depends only on $t$ and is a nondecreasing function of
$t$, $h(u)$ is a $\mathcal{E}$-measurable function which satisfies
the following condition of the {\it Cramer type}: the function
\begin{equation}\label{1.3}
K(\lambda)=\int_E(e^{\lambda h(u)}-1\lambda h(u))q(du)
\end{equation}
is finite for all $\lambda\in R$;

{\bf II}. ($C$-continuity in $X$). For all $u\in E, X\in C$ and
almost all $t\in R_+$ with respect to the Lebesgue measure, the
following holds: id a sequence $X^{(k)}\in D$, $k\ge 1$ is such
that
\[
\lim_{k\to\infty}\sup_{s\le t}|X^{(k)}_s-X_s|=0,
\]
then (as $k\to\infty$)
\[
a(t,X^{(k)})\to a(t,X),\ b(t,X^{(k)})\to b(t,X), \
f(t,X^{(k)},u)\to f(t,X,u).
\]

If we further assume that the differential equation
\begin{equation}\label{1.4}
\dot{Y}_t=a(t,Y), \ Y_0=x
\end{equation}
has a unique solution, then (see Lemma 5.1) we the following
ergodic property: for all $T>0$ as $n\to\infty$
\begin{equation}\label{1.5}
\sup_{t\le T}|X^n_t-Y_t|\xrightarrow{P}0.
\end{equation}
This means that for any open set $A\subset D$ which contains
$(Y_t)_{t\ge 0}$
\[
\lim_nn^{-1}\log P(X^n\in A)=0.
\]

In this paper we study large deviations of $X^n,n\ge 1$, i.e. the
asymptotics (as $n\to\infty$) of $n^{-1}\log P(X^n\in A)$ for sets
$A\in D$ not necessarily containing $Y$.

The problem of large deviations (in $D_{[0,T]}$) for Markov
processes of the form (\ref{1.1}) was considered by Wentzell and
Freidlin \cite{WF} and Wentzell \cite{W}. They assumed that the
functions $a=a(t,x)$, $b=b(t,x)$ and $f=f(t,x,u)$ are bounded and
continuous in $(t,x)$, $x\in R$. Here we prove the large deviation
principle in $D$ ($D_{[0,T]}$ is a consequence) for the case when
$a,b$ and $f$ may depend on on the whole past. We also allow $a,b$
and  $f$ to grow linearly in $X$ and to be only measurable in $t$
(Theorems 2.1 and 2.2). Also we give explicit conditions on $a,b$
and $f$ which ensure the assumptions of our main theorems while
the conditions imposed in \cite{WF} and \cite{W} are not easy to
interpret in terms of the coefficients. We wish to present main
ideas in a most refined form, so we consider the one-dimensional
case. The multidimensional case does not seem to be principally
different. Another distinctive feature of the paper is that we use
a completely new approach based on a counterpart of the Prokhorov
theorem.

For processes without jumps (of the diffusion type), along with
similar lines as ours the result of \cite{WF} and \cite{W} have
been generalized by a number of authors. Stroock \cite{Str},
Azencott \cite{Az}, Baldi, P. and Chaleat-Maurel \cite{Bal},
Narita \cite{Nar} and Friedman \cite{Fried} studied homogeneous
diffusions with unbounded coefficients. The case of
non-homogeneous diffusion with the coefficients depending on the
past in a general manner (but bounded) was considered by Cutland
\cite{Cutl} who invoked the infinitesimal technique. All these
results relate to $C_{[0,T]}$. Micami \cite{Mic} obtained the
result from \cite{W} and \cite{WF} under weaker assumptions.

The case of processes with independent increments was considered
in [13-15].

The paper is organized as follows: in the rest of Section 1 we
outline our approach to obtaining the large deviation principle.
In Section 2 main results (Theorems 2.1. and 2.2) are stated.
Section 3 contains the proof of the fundamental property of
$C$-exponential tightness for $X^n$. In Sections 4 and 5 the
ground is laid for proving upper and lower bounds which are
obtained in Sections 6 and 7. In Section 8 we prove Theorems 2.1
and 2.2. Section 9 gives explicit conditions on the functions
$a,b$ and $f$ which are sufficient for the assumptions Theorems
2.1 and 2.2 to hold.

{\bf 2}. Recall main definitions. We assume that $D$ is supplied
with the Skorokhod-Lindvall metric and is thus a complete
separable metric space. Following Varadhan \cite{Var}, \cite{V} we
say that a sequence $X^n,n\ge 1$ of processes with paths in $D$
obeys the large deviation principle (in $D)$ if

\medskip
0) there exists a function $I=I(\phi),\phi\in D$, with values in
$[0,\infty]$ such that for any $\alpha\ge 0$ the set
$\Phi(\alpha):=\{\phi\in D:I(\phi)\le \alpha\}$ is compact in $D$;

\medskip
1) for any open set $G\subset D$
\[
\mathop{\varliminf}\limits_nn^{-1}\log P(X^n\in G)\ge
-\inf_{\phi\in G}I(\phi);
\]

\medskip
1) for any closed set $F\subset D$
\[
\mathop{\varlimsup}\limits_nn^{-1}\log P(X^n\in F)\le
-\inf_{\phi\in F}I(\phi).
\]

A function $I=I(\phi)$ meeting in (0) is called a rate function
\cite{V} (or good rate function, \cite{DuStr}).

Equivalent to (1) and (2) under (0) are the following conditions
introduced by Wentzell and Freidlin \cite{WF}:

1') for all $\delta>0$ and $\phi\in D$
\[
\mathop{\varliminf}\limits_nn^{-1}\log P(\rho(X^n,\phi)<\delta)\ge
-I(\phi);
\]

2') for all $\delta>0$ and $\alpha\ge 0$
\[
\mathop{\varlimsup}\limits_nn^{-1}\log P(\rho(X^n,\Phi(\alpha)\ge
\delta) \le -\alpha.
\]

The large deviation principle in $D_{[0,T]}$ is defined similarly.

\quad Our approach for proving a large deviation principle relies
on the concept of exponential tightness.

\begin{definition}
(\cite{DuStr}). A sequence $X^n=(X^n_t)_{t\ge 0},n\ge 1$ of
processes with paths in $D$ is said to be exponential tight if for
any $C$ there exists a compact $K_C$ in $D$ such that
\[
\mathop{\varlimsup}\limits_nn^{-1}\log P(X^n\in D\setminus K_C)\le
-C.
\]
\end{definition}

Other names for exponential tightness are ``large deviation
tightness'' \cite{LynSet} and ``strong tightness'' \cite{P1}.

The fundamental importance of this concept was made clear from
Pukhalskii \cite{P1} who proved for large deviations of an analog
of the Prokhorov theorem on the equivalence of weak relative
compactness and tightness for a family of probability measures
(Theorem 1.1 below. We first recall the following definition from
\cite{P1}.

\begin{definition}
A sequence $X^n=(X^n_t)_{t\ge 0},n\ge 1$ of processes with paths
in $D$ is said to obey the partial large deviation principle if
any subsequence $\{n'\}$ of $\{n\}$ contains a further subsequence
$\{n''\}$ such that $X^{n''}$ obeys the obeys the large deviation
principle (with some rate function).
\end{definition}

\begin{theorem}
{\rm(\cite{P1})}. A sequence $X^n=(X^n_t)_{t\ge 0},n\ge 1$ of
processes with paths in $D$ is exponentially tight if and only if
it obeys the partial large deviation principle.
\end{theorem}

\begin{remark}
The theorem hold for any arbitrary Polish space \cite{P1}.
\end{remark}

By Theorem \ref{1.1} the following definition makes sense.

\begin{definition}
An exponentially tight sequence $X^n=(X^n_t)_{t\ge 0},n\ge 1$ is
said to be C-exponentially tight if for any subsequence $X^{n'}$
which obeys the large deviation principle with the rate function
$I'(\phi)$, we have $I'(\phi)=\infty$ for all $\phi\in D\setminus
C$.
\end{definition}

Next we introduce the notion of the local large deviation
principle.

\begin{definition}
A sequence $X^n,n\ge 1$ obeys the local large deviation principle
(in $D$) if for any $\phi=(\phi_t)_{t\ge 0}\in D$ we have
\begin{eqnarray}\label{1.6}
&& \lim_{\delta\to 0}\mathop{\varliminf}\limits_nn^{-1}\log
P(\rho(X^n,\phi)\le \delta)
\notag\\
&=& \lim_{\delta\to 0}\mathop{\varlimsup}\limits_nn^{-1}\log
P(\rho(X^n,\phi)\le \delta) (=-J(\phi))
\end{eqnarray}
(it can be $-\infty=-\infty$).
\end{definition}

It is known \cite{WF} that the large deviation principle implies
the local large deviation principle. The following theorem from
\cite{P1} establishes that under exponential tightness the
converse is also true. \vskip .4in

\begin{theorem}
If a sequence $X^n=(X^n_t)_{t\ge 0},n\ge 1$ of processes with
paths in $D$ is exponentially tight and obeys the local large
deviation principle (in $D$), then the function $J=J(\phi)$ from
(\ref{1.6}) is a rate function, and the sequence $X^n,n\ge 1$
obeys the large deviation principle with the  rate function $I=J$.
\end{theorem}

{\it Proof of Theorem \ref{1.2}.} Assume that an exponential tight
sequence $X^n,n\ge 1$ obeys the local large deviation principle
(\ref{1.6}). By Theorem \ref{1.1}, we can choose a subsequence
$\{\tilde{n}\}$ such that the subsequence $X^{\tilde{n}}$ obeys
the large deviation principle with a rate function $\tilde{I}$. By
\cite{WF} the sequence $X^{\tilde{n}}$ obeys the local large
deviation principle with a rate function $\tilde{I}$. On the other
hand, since
\[
\mathop{\varliminf}\limits_n\le\mathop{\varliminf}\limits_{\tilde{n}}
\le\mathop{\varlimsup}\limits_{\tilde{n}}
\le\mathop{\varlimsup}\limits_n,
\]
we have from (\ref{1.6}) that
\begin{eqnarray}\label{1.7}
&& \lim_{\delta\to 0}\mathop{\varliminf}\limits_{\tilde{n}}
\tilde{n}^{-1}\log P(\rho(X^{\tilde{n}},\phi)\le\delta)
\notag\\
&=&\lim_{\delta\to 0}\mathop{\varlimsup}\limits_{\tilde{n}}
\tilde{n}^{-1}\log P(\rho(X^{\tilde{n}},\phi)\le\delta)=-J(\phi).
\end{eqnarray}
This means that $\tilde{I}=J$ and hence $J$ is a rate function.

Let $G\subset D$ be an open set. Choose a subsequence $\{n'\}$ of
$\{n\}$ such that
\begin{equation}\label{1.8}
\mathop{\varliminf}\limits_nn^{-1}\log P((X^n\in
G)=\lim_{n'}(n')^{-1}\log P(X^{n'}\in G)
\end{equation}
and from $\{n'\}$ by Theorem 1.1 choose a subsequence \{n''\} such
that $X^{n''}$ obeys the the large deviation principle with $I''$.
Then by (1)
\begin{equation}\label{1.9}
\lim_{n''}(n'')^{-1}\log P(X^{n''}\in G)\ge -\inf_{\phi\in
G}I''(\phi).
\end{equation}
Now we take $\{n''\}$ as the subsequence $\{\tilde{n}\}$ above.
Then $I''=J$ and by (\ref{1.8}) and (\ref{1.9})
\[
\mathop{\varliminf}\limits_nn^{-1}\log P(X^n\in G)\ge
-\inf_{\phi\in G}J(\phi).
\]
We have proved (1). (2) is proved in the same way.

To establish (1.6) may be difficult since the Skorokhod metric is
not easy to deal with. However  in our particular case the
following modification of the scheme suggested by Theorem 1.2
applies:

$\alpha$) check that the sequence $X^n,n\ge 1$ is
$C$-exponentially tight,

$\beta$) for all $T>0$ and all $\phi=(\phi_t)_{t\ge 0}\in C$
calculate
\begin{eqnarray}\label{1.10}
J_T(\phi)&=&-\lim_{\delta\to
0}\mathop{\varliminf}\limits_nn^{-1}\log P(\sup_ {t\le
T}|X^n_t-\phi_t|\le\delta)
\notag\\
&=&-\lim_{\delta\to 0}\mathop{\varlimsup}\limits_nn^{-1}\log
P(\sup_ {t\le T}|X^n_t-\phi_t|\le\delta),
\end{eqnarray}

$\gamma$) then
\[
I(\phi)=
  \begin{cases}
    \sup\limits_TJ_T(\phi) & \phi\in C,
    \\
    \infty & \phi\in D\setminus C.
  \end{cases}
\]
This is a consequence of the following theorem.
\begin{theorem}
Let a sequence $X^n,n\ge 1$ be $C$-exponentially tight, and let
(1.10) hold for all $\phi\in C$ and $T>0$. Then $X^n$ obeys the
large deviation principle with the rate function given by
$\gamma${\rm)}.
\end{theorem}

{\it Proof.} By Theorem 1.2 and the definition of $C$-exponential
tightness it suffices to prove that (1.10) implies (1.6) for
$\phi\in C$. First we recall the definition of the
Skorokhod-Lindvall metric in $D$ (see, e.g. \cite{LSMar}).

Let for $k=1,2\ldots,$
\[
g_k(t)=I(t\le k)+(k+1-t)I(k<t\le k+1), \ t\ge 0.
\]
For $X=(X_t)_{t\ge 0}\in D$ and $\phi=(\phi_t)_{t\ge 0}\in D$
define $X^k=(X^k_t)_{0\le t\ge 1}\in D_{[0,1]}$ and
$\phi^k=(\phi^k_t) _{0\le t\ge 1}\in D_{[0,1]}$, $k=1,2,\ldots,$
by
\begin{eqnarray*}
&& X^k_t=X_{\alpha(t)g_k(\alpha(t))}, \ 0\le t\le 1, \ X^k_1=0,
\\
&& \phi^k_t=\phi_{\alpha(t)g_k(\alpha(t))}, \ 0\le t\le 1, \
\phi^k_1=0,
\end{eqnarray*}
where
\begin{equation}\label{1.11}
\alpha(t)=-\log(1-t), \ 0\le t<1.
\end{equation}
Let $d_0$ be the complete metric in $D_{[0,1]}$ introduced by
Prokhorov: if $Y=(Y_t)_{0\le t\le 1}\in D$ and $Z=(Z_t)_{0\le t\le
1}\in D$ then
\begin{equation}\label{1.12}
d_0(Y,Z)=\inf_{\mu\in\frak{M}}\Big\{\sup_{t\le 1}|Y_t-Z_{\mu(t)}|+
\sup_{0\le s<t\le 1}\Big|\log\frac{\mu(t)-\mu(s)}{t-s}\Big|\Big\},
\end{equation}
where $\frak{M}$ is a set of strictly increasing continuous
functions $\mu= (\mu(t))_{0\le t\le 1}$ with $\mu(0)=0$,
$\mu(1)=1$.

The Skorokhod-Lindvall metric is given by
\begin{equation}\label{1.13}
\rho(X,\phi)=\sum_{k=1}^\infty
2^{-k}\frac{\rho_k(X,\phi)}{1+\rho_k(X,\phi)},
\end{equation}
where
\begin{equation}\label{1.14}
\rho_k(X,\phi)=d_0(X^k,\phi^k), \ k=1,2,\ldots
\end{equation}
Since obviously
\[
d_0((Y,Z)\le \sup_{t\le 1}|Y_t-Z_t|,
\]
we have from (1.14) and the definitions of $X^k$ and $\phi^k$ that
for all $k=1,2,\ldots$
\[
\rho_k(X,\phi)\le \sup_{t\le k}|X_t-\phi_t|,
\]
and by (1.13)
\begin{equation}\label{1.15}
\rho(X,\phi)\le \sup_{t\le k+1}|X_t-\phi_t|+\frac{1}{2^k}, \
k=1,2,\ldots
\end{equation}
Next we show that for $\delta\le 1/4$
\begin{equation}\label{1.16}
\{\rho_k(X,\phi)\le\delta\}\subseteq\Big\{\sup_{t\le
k}|X_t-\phi_t|\le 2\delta+ W_k(\phi,4\delta_k)\Big\},
\end{equation}
where
\begin{equation}\label{1.17}
W_k(\phi,\sigma)=\sup_{\mathop{|u-v|\le\sigma}\limits_{0\le u,v\le
k+1}} |\phi_{ug_k(u)}-\phi_{vg_k(v)}|,
\end{equation}
\begin{equation}\label{1.18}
\delta_k=\delta/(1+\alpha^{-1}(k+1)),
\end{equation}
($\alpha^{-1}$ is the inverse of $\alpha$).

First note that if
\[
\sup_{0\le s<t\le 1}\Big|\log\frac{\mu(t)-\mu(s)}{t-s}\Big|\le
\frac{1}{4},
\]
then ([20, ch. 3, \S14])
\begin{equation*}
\frac{1}{2}\sup_{t\le 1}|\mu(t)-t|\le \sup_{0\le s<t\le
1}\Big|\log\frac{\mu(t)-\mu(s)} {t-s}\Big|.
\end{equation*}
Therefore if $\rho_k(X,\phi)\le 1$, then by (1.12) and (1.14)
\begin{equation}\label{1.19}
\rho_k(X,\phi)\ge \inf_{\mu\in\frak{M}}\Big\{\sup_{t\le
1}|X^k_t-\phi^k_{\mu(t)}| +\frac{1}{2}\sup_{t\le
1}|\mu(t)-t|\Big\}.
\end{equation}
So for any $\delta\le \frac{1}{4}$ we can find
$\mu_\delta\in\frak{M}$ such that
\begin{equation}\label{1.20}
\rho_k(X,\phi)\ge \Big\{\sup_{t\le
1}|X^k_t-\phi^k_{\mu_\delta(t)}| +\frac{1}{2}\sup_{t\le
1}|\mu_\delta(t)-t|\Big\}-\delta.
\end{equation}
Thus
\begin{equation}\label{1.21}
\{\rho_k(X,\phi)\le \delta\}\subseteq\Big\{\frac{1}{2}\sup_{t\le
1}|\mu_\delta(t) -t|\le 2\delta\Big\}.
\end{equation}
From (1.11) it is obvious that for $s,t\le \alpha^{-1}(k+1)$
\[
|\alpha(t)-\alpha(s)||t-s|/(1-\alpha^{-1}(k+1))
\]
which easily implies in view of the definition of $X^k$ and
$\phi^k$, and by (1.17), (1.18) and (1.21) that if
$\rho_k(X,\phi)\le\delta$ then
\[
\sup_{t\le 1}|\phi^k_t-\phi^k_{\mu_\delta(t)}|\le
W_k(\phi,4\delta_k).
\]
The latter by (1.20) and by the definition of $X^k$ and $\phi^k$
implies that $\rho_k(X,\phi)\le\delta(\le \frac{1}{4})$ then
\[
\rho_k(X,\phi)\ge \sup_{t\le
1}|X^k_t-\phi^k_t|-W_k(\phi,4\delta_k)-\delta \ge\sup_{t\le
k}|X_t-\phi_t|-W_k(\phi,4\delta_k)-\delta
\]
which gives (1.16).

(1.15) and (1.16) together with (1.10) lead to (1.6) almost
immediately. Indeed from(1.15) and (1.10) we have for all
$\delta>0$ and $k=1,2,\ldots$
\begin{eqnarray*}
&& \mathop{\varliminf}_nn^{-1}\log P(\rho(X^n,\phi)\le
\delta+1/2^k)
\\
&\ge& \mathop{\varliminf}_nn^{-1}\log P\Big(\sup_{t\le
k+1}|X^n_t-\phi_t| \le \delta\Big)\ge -J_{k+1}(\phi).
\end{eqnarray*}
Since by $\gamma)$ $J_{k+1}(\phi)\le J(\phi)$ we obtain
\begin{equation}\label{1.22}
\mathop{\varliminf}\limits_{\delta\to
0}\mathop{\varliminf}\limits_nn^{-1} \log P(\rho(X^n,\phi)\le
\delta)\ge -J(\phi).
\end{equation}
Conversely (1.13) and (1.16) yield for $\delta\le 1/4$ and
$k=1,2,\ldots$
\begin{equation}\label{1.23}
\Big\{\rho(X,\phi)\le\frac{\delta}{1+\delta}\Big/2^k\Big\}\subseteq\{\rho_k(X_k,
\phi)\le\delta\}\subseteq\Big\{\sup_{t\le k}|X_t-\phi_t|\le
2\delta+ W_k(\phi,4\delta_k)\Big\}.
\end{equation}
Since the function $\phi_tg_k(t)$ is continuous, so in view of
(1.17) and (1.18)
$$
W_k(\phi,4\delta_k))\to 0 \  \text{as} \ \delta\to 0,
$$
and then by (1.23) and (1.10),
\begin{eqnarray}\label{1.24}
&& \lim_{\delta\to 0}\mathop{\varlimsup}_nn^{-1} \log
P(\rho(X^n,\phi)\le \delta)
\notag\\
&\le& \lim_{\delta\to 0}\mathop{\varlimsup}_nn^{-1}\log
P\Big(\sup_{t\le k}|X^n_t-\phi_t| \le \delta\Big)\le -J_k(\phi).
\end{eqnarray}
Since by $\gamma$) $J(\phi)=\sup_kJ_K(\phi)$, (1.22) and (1.24)
prove (1.6) for $\phi\in C$.

\medskip
{\it Remark 1.} \ Note that the fact in the theorem has been
noticed by Dawson and G\"artner ([24. Th. 5.1-5.3]) for the
particular problem they studied.

\medskip
{\it Remark 2.} \ Along with the Skorokhod topology one can also
consider the local uniform topology on $D$ (see e.g. [3]) and
study the large deviation principle for this topology, should
$X^n$ remain measurable with respect to the  corresponding
$\sigma$-field, which the case of the solution (1.1). As the local
uniform topology is stronger than the Skorokhod topology, the
large deviation principle for the local uniform topology implies
that for the Skorokhod topology. In fact, the result [4] and [5]
are for the uniform topology on $D_{[0,T]}$. However, as far as
the solution of (1.1) are concerned the corresponding rate
function (for the Skorokhod topology) is infinity at discontinuous
elements of $D$ (Theorems. 2.1 and 2.2 below), and in such a case
the large deviation principles for both topologies are easily seen
to be equivalent (the situation is the same as for weak
convergence, see [20, ch. 3, \S18]). In particular in Theorem 1.3
we have the large deviation principle for the local uniform
topology as well (provided the measurability is preserved).

\section{\bf Cumulant. Legendre-Fenchel transform. Main results}

1. For $\lambda\in R$, $t\in R_+$, $X\in D$ define the cumulant
\begin{equation}\label{2.1}
G(\lambda;t,X)=\lambda a(t,X)+\frac{\lambda^2}{2}b^2(t,X)+
\int_E(e^{\lambda f(t,X,u)}-1-\lambda f(t,X,u))q(du).
\end{equation}
By condition ${\bf I}$, $G(\lambda;t,X)$ is finite and smooth in
$\lambda$ for all $\lambda\in R$ and $X\in D$, and almost all
$t\in R_+$ with respect to the Lebesgue measure, with the first
two derivatives expressed as
\begin{equation}\label{2.2}
\begin{split}
g(\lambda;t,X):=&G'_\lambda(\lambda;t,X)=a(t,X)+\lambda b^2(t,X)+
\int_Ef(t,X,u)(e^{\lambda f(t,X,u)}-1)q(du)
\\
&G''_{\lambda\lambda}(\lambda;t,X)=b^2(t,X)+
\int_Ef^2(t,X,u)e^{\lambda f(t,X,u)}q(du).
\end{split}
\end{equation}
In particular, the cumulant is convex in $\lambda$.

\medskip
{\it Remark.} Throughout the paper we assume that $R_+$ is
supplied by the Lebesgue measure and so in the sequel we omit
specific reference on it.

Define the Legendre-Fenchel transform of $G(\lambda;t,X)$ (cf.,
e.g. [21])
\begin{equation}\label{2.3}
H(y;t,X)=\sup_{\lambda\in R}[\lambda y-G(\lambda;t,X)].
\end{equation}
Since the cumulant is continuous in $\lambda$ ``sup'' in (2.3) may
be taken only over rational $\lambda\in R$ for all $y\in R$, $t\in
R_+$, $X\in D$. Therefore, $H(y;t,X)$ is
$\mathcal{B}(R)\otimes\mathcal{P}(\mathbb{D})$ measurable. Besides
$H(y;t,X)$ is nonnegative since $H(0;t,X)=0$.

\smallskip
2. For our main result we need two additional conditions.

For $\phi=(\phi_t)_{t\ge 0}\in D$ and $T>0$ define
\[
I_T(\phi)=
  \begin{cases}
    \int_0^TH(\dot{\phi}_t;t,\phi)dt, & d\phi_=\dot{\phi}_tdt, t\in[0,T];
    \phi_0=x(\equiv X^n_0), \\
    \infty, & \text{otherwise}.
  \end{cases}
\]
The following condition enables us to calculate ``sup'' in (2.3).

\smallskip
{\bf III} (solvability condition). For all $T>0$ and for all
$\phi\in D$ with $I_T(\phi)<\infty$ there exist
$\delta_{T,\phi}>0$ and a $\mathcal{B}(R)
\otimes\mathcal{P}(\mathbb{D})$-measurable function
$\Lambda_{T,\phi}(y;t,X)$, $y\in R$, $t\in R_+$, $X\in D$, such
that for all $y\in R$ and $X\in D$ with $\sup\limits_{t\le
T}|X_t-\phi_t|\le\delta_{T,\phi}$ we have
\begin{equation}\label{2.4}
y=g(\Lambda_{T,\phi}(y;t,X),t,X)
\end{equation}
for almost all $t\le T$. In additional, $\Lambda_{T,\phi}(y;t,X)$
has the following properties.

\medskip
{\bf i}) (local boundedness). For every $N>0$ there exists $r$
(depending on $N$, $\phi$, $\delta_{T,\phi}$ and $T>0$)), such
that for all $y\in R$ with $|y|\le N$ and $X\in D$ with
$\sup\limits_{t\le T}|X_t-\phi_t|\le\delta_{T,\phi}$ we have
\[
|\lambda_{T,\phi}(y;t,X)|\le r
\]
for almost all $t\le T$.

\smallskip
{\bf ii}) ($C$-continuity). For all $y\in R$,
$\Lambda_{T,\phi}(y;t,X)$ is $C_{[0,T]}$-continuous in $X$ at
$X=\phi\in C$ for almost all $t\le T$, i.e. if $X^{(k)}\in D,k\ge
1$, then the implication holds
\[
\lim_{k\to\infty}\sup_{t\le
T}|X^{(k}-\phi_t|=0\Rightarrow\lim_{k\to\infty}
\Lambda_{T,\phi}(y;t,X^{(k)})=\Lambda_{T,\phi}(y;t,\phi).
\]

\medskip
{\it Remark.} (2.2) implies that under condition {\bf III} ``sup''
in (2.3) is attained at $\lambda=\Lambda_{T,\phi}(y;t,X)$: for all
$y\in R$ and $X\in D$ with $\sup\limits_{t\le
T}|X_t-\phi_t|\le\delta_{T,\phi}$, we have for almost all $t\le T$
\begin{equation}\label{2.5}
H(y;t,X)=\Lambda_{T,\phi}(y;t,X)y-G(\Lambda_{T,\phi}(y;t,X);t,X).
\end{equation}

\medskip
{\bf IV}. For all absolutely continuous functions
$\phi=(\phi_t)_{t\ge 0}$ with $\phi_0=x$ and all $T>0$ the
implication holds
\[
\int_0^TH(\dot{\phi}_t;,t\phi)dt<\infty\Rightarrow
\int_0^T\sup_NH(\dot{\phi}_t;,t\phi^N)dt<\infty,
\]
where
\[
\phi^N_t=x+\int_0^t\dot{\phi}_sI(|\dot{\phi}_s|\le N)ds.
\]

\medskip
{\it Remark.} Conditions {\bf III} and {\bf IV} are basically
conditions of nondegeneracy (see Theorem 9.1 below).

\medskip
4. Now we state main results.

\begin{theorem}
Let conditions {\bf I-IV} hold. Then the sequence $X^n,n\ge 1$, of
semimartingales defined by (1.1) obeys the large deviation
principle in $D$ with the rate function
\[
I(\phi)=
  \begin{cases}
    \int_0^\infty H(\dot{\phi}_t;t,\phi)dt, & \text{$\phi$ is absolutely
    continuous, $\phi_0=x(\equiv X^n_0$},
    \\
    \infty, & \text{otherwise}.
  \end{cases}
\]
\end{theorem}

\medskip
5. Theorem 2.1 does not include the Poisson process. So we
consider separately the following case covers processes of the
Poisson type: for all $T\in R_+$ all $X\in D$ and $u\in E$
\begin{equation}
b(t,X)=0 \ {\rm a.e.},
\end{equation}

\begin{equation}
f(t,X,u)>\epsilon_T \ {\rm a.e. \ on} \ [0,T],
\end{equation}

\begin{equation}
q(E)<\infty.
\end{equation}
\begin{theorem}
Let (2.6)-(2.8) hold and the equation
\begin{equation}
d\psi/dt=a(t,\psi)-\int_Ef(t,\psi,u)q(du), \ \psi_0=x
\end{equation}
has a unique solution. If conditions {\bf I}, {\bf II} and {\bf
IV} are satisfied and condition {\bf III} may fail only for the
function $\phi$ which is a solution of (2.9), then the assertion
of Theorem 2.1 remains true.
\end{theorem}

\medskip
{\it Remark 1.} Since for $X\in D$, $X^*_{t-}=X^*_t$ almost
everywhere in condition {\bf I} $X^*_{t-}$ may be replaced by
$X^*_t$.

\medskip
{\it Remark 2.} Assumption (2.7) in Theorem 2.2 may be replaced by
the assumption $f(t,X,u)\ge q(t,u)>0$ with
\[
\int_0^T\int_Eg^{-1}(t,u)q(du)<\infty.
\]
Analogously one can consider the case when $f$ is negative.

\medskip
{\it Remark 3.} Since $I(\phi)=\infty$ for $\phi\in D\setminus C$
it is easy to deduce applying yje continuous mapping theorem (see
[18, Lemma 2.1.4], [19, Th. 2.2]) that for all $T>0$ the sequence
$((X^n_t)_{0\le t\le T}),n\ge 1$ obeys the large deviation
principle in $D_{[0,T]}$ with the rate function $I_T$ which is
defined above (since $I_T$ depends of the values of $\phi$ up to
$T$, we can as well consider is as a function on $D_{[0,T]}$.

\medskip
{\it Remark 4.} All the results are retained for the (local)
uniform topology in $D(D_{[0,T]}$) (see the remark at the end of
Section 1).

\section{\bf Exponential tightness}

1. In this section we prove that under condition {\bf I} the
sequence $X^n,n\ge 1$, defined by (1.1), is $C$-exponentially
tight. For this, we use the next theorem on $C$-exponential
tightness of a sequence of adapted processes.

\begin{theorem}
Let $X^n=(X^n_t)_{t\ge 0},n\ge 1$ be a sequence of processes with
paths in $D$. Each $X^n$ is defined on a stochastic basis
$(\Omega,\mathcal{F}, \mathbb{F}^n,P)$. Assume that the following
is satisfied: for all $L>0$ and $\eta
>0$
\begin{eqnarray*}
&& i)\quad \lim_{c\to\infty}\mathop{\varlimsup}\limits_nn^{-1}\log
P(X^{n^*}_L\ge c) =-\infty,
\\
&& ii)\quad \lim_{\delta\to
0}\mathop{\varlimsup}\limits_n\sup_{\tau\in T_L(\mathbb{F}^n)}
n^{-1}\log P(\sup_{t\le \delta}|X^n_{\tau+t}-X^n_\tau|\ge
\eta)=-\infty.
\end{eqnarray*}

Then the sequence $X^n,n\ge 1$ is $C$-exponentially tight.
\end{theorem}

\medskip
{\it Proof.} By Theorem 4.4 [19], i) and ii) imply that the
sequence $X^n,n\ge 1$ is exponentially tight. By Theorem 4.6 [19]
it is proved that ii) also implies $C$-exponential tightness. We
reproduce the prof here in more details.

Let $\phi\in D\setminus C$, We can find $k\in\{2,3,\ldots\}$ such
that there exists $s_0\in (0,k-1]$ with $\triangle\phi_{s_0}\neq
0$. Let us show that for $Y\in D$ we have if $\delta\le
(1\wedge(\alpha^{-1}(k)-\alpha^{-1}(k-1))/4$ (we use the notation
of the proof of Theorem 1.3)
\begin{equation}
\{\rho_k(Y,\phi)\le\delta\}\subseteq\Big\{\sup_{|s_0-t|\le
4\delta(k)}|Y_{s_0}- Y_t|\ge |\triangle\phi_{s_0}|-2\delta\Big\},
\end{equation}
where $\delta(k)=\delta/(1-\alpha^{-1}(k)).$

As in the proof of Theorem 1.3, we can choose
$\mu_\delta=(\mu_\delta(t))_ {0\le t\le 1}\in\frak{M}$. so that
\[
\rho_k(Y,\phi)\ge \sup_{t\le 1}|Y^k_{\mu_\delta(t)}-\phi^k_t|+
\frac{1}{2}\sup_{t\le 1}|\mu_\delta(t)-t|-\delta.
\]
Then if $\rho_k(Y,\phi)\le \delta$ we have
\[
\frac{1}{2}\sup_{t\le 1}|\mu_\delta(t)-t|\le 2\delta
\]
and hence, taking $t_0=\alpha^{-1}(s_0)$
\begin{equation}
|\triangle Y^k_{\mu_\delta(t_0)}|\le 2\sup_{\mathop{|t-s|\le
4\delta}\limits _{0\le s,t\le 1}}|Y^k_{t_0}-Y^k_t|.
\end{equation}
As by the triangular inequality
\[
|\triangle\phi^k_{t_0}|\le|\triangle
Y^k_{\mu_\delta(t_0)}-\phi^k_{t_0}| +|\triangle
Y^k_{\mu_\delta(t_0)}|\le 2\sup_{t\le 1}|\triangle
Y^k_{\mu_\delta(t_0)}-\phi^k_t|+ |\triangle
Y^k_{\mu_\delta(t_0)}|,
\]
so we deduce from (3.2) that if $\rho_k(Y,\phi)\le \delta$, then
\[
\sup_{\mathop{|t_0-t|\le 4\delta}\limits_{0\le t\le
1}}|Y^k_{t_0}-Y^k_t|\ge |\triangle\phi^k_{t_0}|/2-\sup_{t\le
1}|Y^k_{\mu_\delta(t)}-\phi^k_t|,
\]
and therefore
\[
\sup_{\mathop{|t_0-t|\le 4\delta}\limits_{0\le t\le
1}}|Y^k_{t_0}-Y^k_t|\ge
|\triangle\phi^k_{t_0}|/2-(\rho_k(Y,\phi)+\delta).
\]
Now $\alpha(t_0)\le k-1$ and $4\delta\le
\alpha^{-1}(k)-\alpha^{-1}(k-1)$. So by the definition of $Y^k$
(see the proof of Theorem 3.1)
\[
\sup_{\mathop{|t_0-t|\le 4\delta}\limits_{0\le t\le
1}}|Y^k_{t_0}-Y^k_t|= \sup_{\mathop{|t_0-t|\le
4\delta}\limits_{0\le t\le 1}}|Y^k_{\alpha(t_0)}-
Y^k_{\alpha(t)}|.
\]
Noticing that
\[
|\alpha(t)-\alpha(t_0)|\le 4\delta/(1-\alpha^{-1}(k))
\]
for $|t-t_0|\le 4\delta$, we obtain (3.1).

Obviously since $s_0<k$
\[
P\Big(\sup_{|s_0-t|\le 4\delta(k)}|X^n_{s_0}-X^n_t|\ge \eta\Big)
\le 2\sup_{\tau\le T_k(\mathbb{F}^n)}P\Big(\sup_{s\le
4\delta(k)}|X^n_{\tau+s} -X^n_\tau|\ge \eta\Big).
\]
Therefore, using (3.1) and ii) we conclude that
\begin{eqnarray*}
&& \lim_{\delta\to 0}\mathop{\varlimsup}\limits_nn^{-1}\log
P(\rho_k(X^n,\phi)\le \delta)
\\
&=& \lim_{\delta\to 0}\mathop{\varlimsup}\limits_nn^{-1}\log\Big\{
2\sup_{\tau\le T_k(\mathbb{F}^n)}P\Big(\sup_{s\le
4\delta(k)}|X^n_{\tau+s} -X^n_\tau|\ge
|\phi_{s_0}|/2-2\delta\Big)\Big\}=-\infty.
\end{eqnarray*}
Since $\rho(X^n,\phi)\ge
2^{-k}\rho_k(X^n,\phi)/(1+\rho_k(X^n,\phi))$ it then follows that
\[
\lim_{\delta\to 0}\mathop{\varlimsup}\limits_nn^{-1}\log
P(\rho(X^n,\phi)\le \delta)=-\infty, \ \phi\in D\setminus C.
\]
Now assume that a subsequence $X^n$ obeys the large deviation
principle wit a rate function $I'$. Then as it obeys the local
large deviation principle, by the above we have $I'(\phi)=\infty,
\phi\in D\setminus C$.

\smallskip
2. Now we state the main result on $C$-exponential tightness of
the sequence of solutions of (1.1).

\begin{theorem}
Under condition {\bf I} (linear growth) the sequence $X^n,n\ge 1$
of solutions of (1.1) is $C$-exponentially tight in $D$.
\end{theorem}

\smallskip
3. The proof is based on a number of lemmas.

\begin{lemma}
Let $Y=(T_t)_{t\ge 0}$ be a semimartingale. Denote by $(B,C,\nu)$
the triplet of predictable characteristic of $Y$ {\rm([2])}, and
assume that $\nu$ satisfies the following analogue of the Cramer
condition
\begin{equation}
\int_0^t\int_{R_0}[e^{\lambda x}-1-\lambda x]\nu(ds,dx)<\infty, \
P-a.s., \ t>0, \ \lambda\in R.
\end{equation}
Assume that for $T>0$ there exists a convex function
$H=H(\lambda),\lambda\in R$ with $H(0)=0$ and such that for all
$\lambda\in R$ and $t\le T$
\[
\lambda\widetilde{B}_t+\lambda^2C_t/2+\int_0^t\int_{R_0}[e^{\lambda
x}-1- \lambda x]\nu(ds,dx)\le H(\lambda\xi), \ P-a.s.,
\]
where $\widetilde{B}_t=B_t+\int_0^t\int_{R_0}xI(|x|>1)\nu(ds,dx)$
and $\xi$ is a nonnegative random variable defined on the same
probability space as $Y$.

Then for all $c>0$ and $\eta>0$
\[
P(Y^*_T\ge\eta)\le P(\xi>c)+\exp\Big\{-\sup_{\lambda\in
R}[\lambda\eta- TH(\lambda\xi)]\Big\}.
\]
\end{lemma}

The proof is  analogous to the proof of Theorem 4.13.2 in [2].

\begin{lemma}
Assume that $X^n=(X^n_t)_{t\ge 0}$ for each $n$ is a special
semimartingale on a stochastic basis
$(\Omega,\mathcal{F},\mathbb{F}^n,P)$ with the decomposition
\begin{equation}
X^n_t=x+A^n_t+M^n_t, \ x\in R.
\end{equation}
Assume that the predictable process of locally bounded variation
$A^n=(A^n_t)_ {t\ge 0}$ and the local martingale
$M^n=(M^n_t)_{t\ge 0}$ have the following properties:

1) for all $T>0$ there exists $a_0=a_0(T)$ such that for all $t\le
T$
\begin{equation}
Var(A^n)_t\le a_0\int_0^t(1+X^{n*}_s)ds;
\end{equation}

2) $(M^n)^{2n}((M^n_t)^{2n})_{t\ge 0}$ is a special semimartingale
with the decomposition
\begin{equation}
(M^n_t)^{2n}=V^n_t+L^n_t,
\end{equation}
where $L^n=(L^n_t)_{t\ge 0}$ is a local martingale, and
$V^n=(V^n_t)_{t\ge 0}$ is an increasing predictable process such
that for all $T>0$ there exists $a_1=a_1(T)$ and $a_2=a_2(T)$ for
which
\begin{equation}
V^n_t\le(a_1+a_2n)\int_0^t\{1+(M^{n*}_s)^{2n}\}ds
\end{equation}
if $t\le T$. Then for all $T>0$
$$
\mathop{\varlimsup}\limits_nn^{-1}\log E(X^{n*}_T)^{2n}<\alpha,
$$
where $\alpha$ depends only on $T$, $a_0$, $a_1$, $a_2$.
\end{lemma}

{\it Proof.} If $s\le t\le T$, then by (3.4) and (3.5)
\[
X^{n*}_s\le |x|++a_0\int_0^s(1+X^{n*}_u)du+M^{n*}_s.
\]
This and the Gronwall inequality imply that
\begin{equation}
X^{n*}_t\le (|x|+a_0t+M^{n*}_t)\exp(a_0t)
\end{equation}
and hence, by Jensen inequality
\[
(X^{n*}_t)^{2n}=3^{2n-1}(|x|^{2n}+(a_0t)^{2n}+((M^{n*}_t)^{2n})\exp(2na_0t),
\ t\le T.
\]
It then follows that for fixed $T$ we can find a constant $C$
which depends on $x,a_0$ and $T$, and such that
\[
(X^{n*}_T)^{2n}\le 1/2(C)^n\{1+(M^{n*}_T)^{2n}\}.
\]
Hereafter
\begin{equation}
E(X^{n*}_T)^{2n}\le 1/2(C)^n\{1+E(M^{n*}_T)^{2n}\}\le
1/2(C)^n\{1+\vee E(M^{n*}_T)^{2n}\}.
\end{equation}
Now we estimate $E(M^{n*}_T)^{2n}$. Assume that $t\le T.$ By the
version of the Doob inequality in [2] (Th. 19.2), decomposition
(3.6) and inequality (3.7) we have
\[
E(M^{n*}_t)^{2n}\le(2n/(2n-1))^{2n}V^n_t\le
4(a_1+a_2n)\int_0^t\{1+e(M^{2*}_s)^{2n}\}ds.
\]
Therefore, if $E(M^{n*}_T)^{2n}<\infty$, then the Gronwall
inequality yields
\[
E(M^{n*}_T)^{2n}\le (a_1+a_2n)T\exp\{4(a_1+a_2n)T\}
\]
(in general case, (3.6) implies the local integrability of
$((M^{n*}_t)^{2n}$ and the above inequality follows by
localization).

The latter and (3.9) imply that
\[
\mathop{\varlimsup}\limits_n1/n\log E(X^{n*}_T)^{2n}\le\le log
C+4a_2T \ (=\alpha).
\]
\smallskip

4. We prove Theorem 3.2 by verifying conditions i) and ii) of
Theorem 3.1.

To check i) we use Stroock's idea ([6]), Lemma 4.12, ch. 4). By
the Chebyshev inequality
\[
P(X^{n+}_L>c)\le c^{-2n}E(X^{n*}_L)^{2n}.
\]
Hence i) holds if
\begin{equation}
\mathop{\varlimsup}\limits_n1/n\log E(X^{n*}_L)^{2n}<\infty.
\end{equation}
To prove (3.10), we apply Lemma 3.2. By (1.1) we have
decomposition (1.4) in which
\begin{equation}
A^n_t=\int_0^ta(s,X^n)ds,
\end{equation}

\begin{equation}
M^n_t=1/\sqrt{n}\int_0^tb(s,X^n)dW^n_s+1/n\int_0^t\int_Ef(s,X^n,u)[p^n-q^n]
(ds,du).
\end{equation}

Conditions (3.5) of Lemma 3.2 in view of condition {\bf I} (of
linear growth).

Now, we verify (3.6) and (3.7). By the It\^o formula (cf,, e.g.
[2], [3]),
\begin{equation}
(M^n_t)^{2n}=2n\int_0^t\int_0^t(M^n_{s-})^{2n-1}dM^n_s+(2n-1)\int_0^t
(M^n_s)^{2n-2}b^2(s,X^n)ds+U^n_t,
\end{equation}
where
\begin{equation}
U^n_t=\sum_{s\le
t}\{(M^n_s)^{2n}-(M^n_{s-})^{2n}-2n(M^n_{s-})^{2n-1}\triangle
M^n_s\}
\end{equation}
is an increasing process. If the process $U^n=(U^n_t)_{t\ge 0}$ is
locally integrable, then we have (3.6) with
\begin{equation}
V^n_t=(2n-1)\int_0^t(M^n_s)^{2n-2}b^2(s,X^n)ds+\widetilde{U}^n_t
\end{equation}
and
\[
L^n_t=2n\int_0^t(M^n_{s-})^{2n-1}dM^n_s+U^n_t-\widetilde{U}^n_t,
\]
where $\widetilde{U}^n=(\widetilde{U}^n_t)_{t\ge 0}$ is the
compensator of $U^n$. Let us show that $U^n$ is locally
integrable. By (1.1), $\triangle X^N_s\equiv \triangle M^n_s$
where in view of (1.1) and (1.2)
\[
\triangle X^n_s=1/n\int_Ef(s,X^n,u)p^n(\{s\},,du).
\]
Thus if we denote
\[
H^n(s,u)=\{M^n_{s-}+1/nf(s,X^n,u)\}^{2n}-(M^n_{s-})^{2n}-2(M^n_{s-})^{2n-1}
f(s,X^n,u),
\]
then
\[
U^n_t=\int_0^tH^n(s,u)p^n(ds,du).
\]
Since $H^n(s,u)\ge 0$ and $H^n(s,u)$ is
$\widetilde{\mathcal{P}}$-measurable, and $q^n(ds,du)$ is the
compensator of $p^n(ds,du)$, the process $U^n$ is locally
integrable and the compensator
$\widetilde{U}^n=(\widetilde{U}^n_t)_{t\ge 0}$ of $U^n$ is
\begin{equation}
\widetilde{U}^n_t=\int_0^t\int_EH^n(s,u)q^n(ds,du),
\end{equation}
provided the integral in (3.16) is finite $P$-a.s. for all $t>0$.
Indeed, for any stopping time $\tau$
\begin{equation}
EU^n_\tau=E\widetilde{U}^n_\tau
\end{equation}
(cf., e.g. Th. 3.2.1 in [2]). Furthermore, $\widetilde{U}^n$ being
an increasing finite-valued predictable process, is locally
bounded (Lemma 1.6.1 in [2]). It follows then that (3.17) implies
both the local integrability of $U^n$ and the fact that
$\widetilde{U}^n$ is the compensator of $U^n$.

Now we prove that that integral in (3.16) is finite $P$-a.s. To
this end denote
\begin{equation}
a=1\vee M^{n*}_{s-},\quad b=|1/nf(s,X^n,u)|.
\end{equation}
Applying the mean value theorem to $H^n(s,u)$ we have
\begin{eqnarray}
H^n(s,u)\le
n(2n-1)(a+b)^{2n-2}b^2&=&n(2n-1)a^{2n-2}(1+b/a)^{2n-2}b^2
\notag\\
&\le&n(2n-1)a^{2n-2}\exp\{(2n-2)b/a\}b^2.
\end{eqnarray}
Since $X^n$ admits a representation (3.4) with $A^n$ from (3.11)
and $M^n$ from (3.12), and $A^n$ by condition {\bf I} satisfies
(3.5), it follows from the proof of Lemma 3.2 that for $X^{n*}_t$
an inequality of the type (3.8) holds. Hence, there exists $c>1$
such that for $s\le T$
\begin{equation}
X^{n*}_{s-}\le c(1\vee M^{n*}_{s-}).
\end{equation}
Therefore (see (3.18)), $a\ge 1\vee (c^{-1}X^{n*}_{s-}) \ (=a')$.
So by (3.19), we have for $n\ge 2$
\begin{equation}
H^n(s,u)\le
n(2n-1)(2n-2)^{-2}a^{2n}\exp\{(2n-2)b/a'\}[(2n-2)b/a']^2.
\end{equation}
Now, we estimate $(2n-2)b/a'$. By condition {\bf I}
\begin{eqnarray*}
(2n-2)b/a'&=&(2n-2)n^{-1}|f(s,X^n,u)|/[1\vee(c^{-1}X^{n*}_{s-})]
\\
&\le&(2n-2)n^{-1}l_sh(u)(1+X^{n*}_{s-})/[1\vee(c^{-1}X^{n*}_{s-})].
\end{eqnarray*}
Since (recall that $c>1$)
\[
1\vee (c^{-1}X^{n*}_{s-})\ge 1/2(1+c^{-1}X^{n*}_{s-})\ge(2c)^{-1}
(1+X^{n*}_{s-})
\]
we have
\[
(2n-2)b/a'\le 4cl_sh(u)\le 4cl_Th(u) \ (=r_Th(u)).
\]
So (3.2) implies that
\[
H^n(s,u)\le 2a^{2n}(r_Th(u))^2\exp(r_Th(u)), \ n\ge 2.
\]
The latter allows us to obtain an estimate for the integral in
(3.16). We have (see (3.18)) and the definition of $a$)
\[
\widetilde{U}^n_t\le 2\int_0^t\int_E(1\vee
M^{n*}_{s-})^{2n}(r_Th(u))^2 \exp(r_Th(u))q^n(ds,du).
\]
Choose $\beta>0$ such that the inequality $x2e^x\le e^{\beta
x}-1-\beta x$ holds for all $x>0$. Then by (1.3)
\[
\gamma=\int_E(r_Th(u))^2\exp(r_Th(u))q(du)<\infty.
\]
Then, since $q^n((ds,du)$ equals $ndsq(du)$ (see (1.2)), so
\begin{equation}
\widetilde{U}^n_t\le 2\gamma nt(1+M^{n*}_{s-})^{2n}<\infty, \
P-\text{a.s.}
\end{equation}

So we have established that $U^n$ is locally integrable, and as a
sequence, (3.6) holds with $V^n$ from (3.15).

Now we prove (3.7). To this end note that in view of (3.22)), the
process $V^n$ from (3.15) can be estimated as
\begin{equation}
V^n_t\le (2n-1)\int_0^t(M^{n*}_s)^{2n-2}b^2(s,X^n)ds+2\gamma
n\int_0^t (1+M^{n*}_s)^{2n}ds.
\end{equation}
By condition {\bf I}, $b^2(s,X^n)\le l_T(1+X^{n*}_s)^2$, so using
(3.20), we deduce that $V^n_t$ satisfies (3.7) for some $a_1$ and
$s_2$ which depend only on $T$.

Thus, all the conditions of Lemma 3.2 are met, and hence, (3.10)
holds. This ends the proof of condition i) of Theorem 3.1.

Check ii). Denote
\[
\quad \ \ A_{n,\delta,\tau}=\Big\{\sup_{t\le\delta}|X^n_{\tau+t}-
X^n_\tau|>\eta\Big\}
\]
and
\begin{eqnarray*}
&&
A^+_{n,\delta,\tau}=\Big\{\sup_{t\le\delta}(X^n_{\tau+t}-X^n_\tau)>\eta\Big\}
\\
&&
A^-_{n,\delta,\tau}=\Big\{\sup_{t\le\delta}(X^n_\tau-X^n_{\tau+t})>\eta\Big\}.
\end{eqnarray*}
Since
\[
P(A_{n,\delta,\tau})=P(A^+_{n,\delta,\tau}\cup
A^-_{n,\delta,\tau})) \le 2 \{P(A^+_{n,\delta,\tau})\vee
P(A^-_{n,\delta,\tau})\},
\]
it suffices to show that
\begin{equation}
\lim_{\delta\to 0}\mathop{\varlimsup}\limits_n1/n\log\sup_{\tau\in
T_L(\mathbb{F} ^n}P(A^\pm_{n,\delta,\tau})=-\infty.
\end{equation}
We are going to establish (3.24) for $A^+_{n,\delta,\tau}$. We
shall apply Lemma 3.1. Denote $Y^n_t=X^n_{\tau+t}-X^n_\tau$. The
process $Y^n_t=(Y^n_t)_ {t\ge 0}$ is a semimartingale with respect
to the filtration $(\mathcal{F}^n _{\tau+t})_{t\ge 0}$, and by
(1.1) and (1.2), its triplet of predictable characteristics
$(B^n,C^n,\nu^n)$ is as follows
\begin{eqnarray*}
&&B^n_t=\int_\tau^{\tau+t}a(s,X^n)ds-\int_\tau^{\tau+t}\int_Ef(s,X^n,u)
I(|f(s,X^n,u)|>n)q(du)du,
\\
&&\Big(\widetilde{B}^n_t=\int_\tau^{\tau+t}a(s,X^n)ds\Big)
\\
&&C^n_t=n^{-1}\int_\tau^{\tau+t}b^2(s,X^n)ds
\\
&&\nu^n((0,t],\Gamma)=n\int_\tau^{\tau+t}\int_RI(f(s,X^n,u)/n\in\Gamma)q(du)ds,
\ \Gamma\in \mathcal{B}(R_0).
\end{eqnarray*}

Then
\begin{eqnarray*}
&&
\lambda\widetilde{B}^n_t+\lambda^2/2C^n_t+\int_0^t\int_E(e^{\lambda
x}-1 -\lambda x)\nu^n(ds,dx)
\\
&=&\lambda\int_\tau^{\tau+t}a(s,X^n)ds+\lambda^2/2n\int_\tau^{\tau+t}b^2(s,X^n)ds
\\
&&+n\int_\tau^{\tau+t}\int_E(e^{\lambda f(s,X^n,u)}-1-\lambda
f(s,X^n,u)) q(du)ds.
\end{eqnarray*}
Condition {\bf I} implies that $Y^n$ satisfies the assumptions of
Lemma 3.1 for $T=\delta\le 1$ with the function
\[
H^(\lambda)=|\lambda|+\lambda^2/2n+n\int_E(e^{|\lambda|h(u)/n}-1-|\lambda|h(u)/n)q(du)
\]
and the random variable $\xi=l_{L+1}(1+X^{n*}_{L+1}).$

Choose arbitrary $c>0$. We have
\[
\sup_{\lambda\in R}[\lambda\eta-\delta H^n(\lambda
c)]=\sup_{\lambda>0}[\lambda\eta- \delta H^n(\lambda
c)]=n\mathcal{L}(\delta),
\]
where
\[
\mathcal{L}(\delta)=\sup_{\mu>0}\Big[\mu(\eta-\delta
c)-\delta\mu^2c^2/2-\delta \int_E(e^{c\mu h(u)}-1-c\mu
h(u))q(du)\Big].
\]
It is easy to see that
\begin{equation}
\lim_{\delta\to 0}\mathcal{L}(\delta)=\infty.
\end{equation}
Hence by lemma 3.1 we have for $\delta<1$
\begin{eqnarray*}
P\Big\{\sup_{t\le \delta}(X^n_{\tau+t}-X^n_\tau)>\eta)&\le&
P\{l_{L+1}(1+X^N*_{L+1})>c\}+e^{-n\mathcal{L}(\delta)}
\\
&\le&2[P\{X^{n*}>c/l_{L+!}-1\}\vee e^{-n\mathcal{L}(\delta)}.
\end{eqnarray*}
This in view of (3.25) yields
\begin{eqnarray*}
&& \mathop{\varlimsup}\limits_{\delta\to
0}\mathop{\varlimsup}\limits_n1/n\log P\Big\{\sup_{t\le
\delta}(X^n_{\tau+t}-X^n_\tau)>\eta\Big\}
\\
&\le&\mathop{\varlimsup}\limits_n1/n\log
P\{X^{n*}_{L+1}>c/l_{L+1}-1\}.
\end{eqnarray*}
And consequently, (3.24) holds for $A^+_{n,\delta,\tau}$, as we
proved earlier
\[
\lim_{c\to\infty}\mathop{\varlimsup}\limits_n1/n\log
P\big\{X^{n*}_{L+1}>c/l_{L+1}-1\big\}= -\infty.
\]

(3.24) for $A^-_{n,\delta,\tau}$ is established similarly.

\section{\bf Multiplicative decomposition. Change of measure}

1. Let $X^n=(X^n_t)_{t\ge 0}$ be a semimartingale defined in
(1.1). In what follows, one of the main actors is the process
\begin{equation}
U^n_t(\lambda)=\exp\Big(n\int_0^t\lambda(s)xX^n_s\Big),
\end{equation}
where $\overline{\lambda}=(\lambda(t))_{t\ge 0}, \ \lambda(t)|\le
r$ for almost all $t$, is a predictable bounded function.

Applying the It\^o formula, it is easy to see that the process
$Z^n(\lambda)=(Z^n_t(\lambda))_ {t\ge 0}$ with
\begin{equation}
Z^n_t(\overline{\lambda})=
\exp\Big(n\Big[\int_0^t\lambda(s)dX^n_s-\int_0^tG(\lambda(s);s,X^n)ds\Big]\Big),
\end{equation}
where $G(\lambda;t,X)$ is the cumulant (see (2.3)), is a positive
local martingale. So we have the following multiplicative
decomposition
\begin{equation}
U^n_t(\overline{\lambda})=Z^n_t(\overline{\lambda})
\exp\Big(n\int_0^tG(\lambda(s);s,X^n)ds\Big]\Big)
\end{equation}
(cf., [2], ch. 2).

Let $\phi=(\phi_t)_{t\ge 0}\in D$, $\phi_0=x \ (\equiv X^n_0)$ and
$T>0$. For $\gamma>0$ define a stopping time by
\begin{equation}
\tau=\inf\Big(t\le T: \sup_{s\le t}|X^n_s-\phi_s|>\gamma\Big), \
(\inf\{\varnothing\}=T),
\end{equation}
and introduce the stopped process $Z^{n,\tau}(\overline{\lambda})=
(Z^{n,\tau}_t(\overline{\lambda}))_{t\ge 0}$ with
$Z^{n,\tau}_t(\overline{\lambda})= Z^n_{t\wedge
\tau}(\overline{\lambda})$. Obviously
$Z^{n,\tau}_\infty(\overline{\lambda})=Z^{n,\tau}_\tau(\overline{\lambda})$.

\begin{lemma}
Under condition {\bf I}, $Z^{n,\tau}(\overline{\lambda})$ is a
square integrable martingale with
$$
EZ^{n,\tau}_\infty(\overline{\lambda})=1.
$$
\end{lemma}

\medskip
{\it Proof.} Show that
\begin{equation}
E(Z^{n,\tau}_\infty(\overline{\lambda})^2<\infty.
\end{equation}
Using (4.3) (with $2\overline{\lambda}$ and (4.3) we have
\begin{eqnarray*}
&& E(Z^{n,\tau}_\infty(\overline{\lambda}))^2
\\
&=& E(Z^{n,\tau}_\tau(\overline{\lambda}))^2
=E\Big[U^n_\tau(2\overline{\lambda})\exp\Big(-2n\int_0^\tau
G(\overline{\lambda} (s);,s,X^n)ds\Big)\Big]
\\
&=&E\Big\{Z^n_\tau(2\overline{\lambda})\exp\Big[n\int_0^\tau
G(2\overline{\lambda} (s);,s,X^n)ds-2n\int_0^\tau
G(\overline{\lambda}(s);s,X^n)ds\Big]\Big\}.
\end{eqnarray*}
The process
$Z^n(2\overline{\lambda})=(Z^n_t(2\overline{\lambda}))_{t\ge 0}$
is a positive local martingale and hence, it is a positive
supermartingale ([2]), Problem 1.4.4). So
$EZ^n_\tau(2\overline{\lambda})\le 1$, and (4.5) holds, provided
$|G(\alpha;s,X)|$ is uniformly bounded in $s$ for almost all $s\le
\tau$ for any $\alpha>0$. The latter is a consequence of condition
{\bf I} since
\[
|G(\alpha;s,X^n)|\le \alpha l_T(1+\gamma+\phi^*_T)^2+K(\alpha
l_T(1+\gamma+ \phi^*_T)), \ s<\tau.
\]
(4.5) implies that $Z^{n,\tau}(\overline{\lambda})$ is a uniformly
integrable martingale and hence
$$
EZ^{n,\tau}_\infty(\overline{\lambda})=EZ^{n,\tau}_0(\overline{\lambda})=1.
$$

\smallskip
2. Assume that conditions {\bf I} and {\bf III} hold,
$\phi=(\phi_t)_{t\ge 0}$ is absolutely continuous function with
$\phi_0=x \ (\equiv X^n_0)$ and its derivative $\dot{\phi}_t$ is
bounded almost everywhere: $|\dot{\phi}_t\le N$.

For $T>0$ define the stopping time $\tau$ by (4.4) with
$\gamma=\delta_{T,\phi}$ from condition {\bf III}. We set
\begin{equation}
\overline{\lambda}(t)=\Lambda_{T,\phi}(\dot{\phi}_t;t,X^n)I(t\le\tau),
\end{equation}
where $\Lambda_{T,\phi}(y;t,X)$ is defined in condition {\bf III}.
By {\bf III(i)} $\overline{\lambda}=(\overline{\lambda}(t),t\ge
0)$ from (4.6) ia predictable and bounded (almost everywhere in
$t$). Let $Z^{n,\tau}(\overline{\lambda})$ be defined by (4.2)
Then by {\bf III(i)} and Lemma 4.1 $EZ^{n,\tau}_\infty=1$ (in what
follows we fix $\overline{\lambda}$ from (4.6) and omit dependence
on it).

On the measurable space $(\Omega,\mathcal{F})$ define a new
probability measure $Q^{n,\phi}$ by
\begin{equation}
dQ^{n,\phi}=Z^{n,\tau}_\infty dP.
\end{equation}
The definition of $Z^{n,\tau}$ implies that $Z^{n,\tau}_\infty>0$,
$P$-a.s. So the measures $P$ and $Q^{n,\phi}$ are equivalent and
\begin{equation}
dP=(Z^{n,\tau}_\infty)^{-1}dQ^{n,\phi}.
\end{equation}
Since an absolutely continuous change measure preserves the
semimartingale property`[2], $X^n$ is a semimartingale on the
basis $(\Omega,\mathcal{F}, \mathbb{F}^n,Q^{n\phi})$. The
structure of $X^n$ is given in the next theorem.
\begin{theorem}
Let conditions {\bf I} and {\bf III(i)} hold. Then $X^n$ on
$(\Omega,\mathcal{F}, \mathbb{F}^n,Q^{n\phi})$ has the
decomposition
\begin{equation}
X^n_t=x+\int_0^tI(s>\tau)a(s,X^n)ds+M^{n,\phi}_t,
\end{equation}
where $M^{n,\phi}=(M^{n,\phi}_t)_{t\ge 0}$ is a local square
integrable martingale with the predictable quadratic variation
process
\begin{eqnarray}\
\langle
M^{n,\phi}\rangle_t&=&n^{-1}\int_0^tb^2(s,X^n)ds+n^{-1}\int_0^t\int_E
f^2(s,X^n,u)
\notag\\
&&\times\exp[\Lambda_{T,\phi}\dot{\phi}_s,s,X^n)I(s\le\tau)f(s,X^n,u)]q(du)ds.
\end{eqnarray}
\end{theorem}

\medskip
{\it Proof.} By Lemma 4.1, $Z^{n,\tau}$ is a square integrable
martingale. From the definition of $Z^n_t$, we have by the It\^o
formula
\begin{equation}
Z^{n,\tau}_t=1+\int_0^tI(s\le\tau)Z^{n,\tau}_{s-}dL^n_s,
\end{equation}
where $L^n=(L^n_t)_{t\ge 0}$ is a local square integrable
martingale with
\begin{eqnarray}
L^n_t&=&\sqrt{n}\int_0^{t\wedge\tau}\Lambda_{T,\phi}(\dot{\phi}_s,s,X^n)
b(s,X^n)dW^n_s
\notag\\
&&
+\int_0^{t\wedge\tau}\int_E\{\exp[\Lambda_{T,\phi}(\dot{\phi}_s,s,X^n)f(s,X^,u)]
-1\}[p^n-q^n]9ds,du).
\end{eqnarray}
By (1.1), $X^n_t$ has the following decomposition (with respect to
$(\mathbb{F}^n,P)$:
\begin{equation}
X^n_t=x+\int_0^ta(s,X^n)ds+M^n_s,
\end{equation}
where the local square integrable martingale $M^n=(M^n_t)_{t\ge
0}$ is of the form:
\begin{equation}
M^n_t=\frac{1}{\sqrt{n}}\int_0^tb(s,X^n)dW^n_s+\frac{1}{n}\int_0^t\int_E
f(s,X^n,u)[p^n-q^n]9ds,du).
\end{equation}
It follows from (4.11), (4.12) and (4.14) that the predictable
quadratic covariance process $\langle Z^{n,\tau},M^n\rangle$ is
given by
\begin{eqnarray*}
\langle
Z^{n,\tau},M^n\rangle_t&=&\int_0^{t\wedge\tau}Z^{n,tau}_{s-}\Big[
\Lambda_{T,\phi}(\dot{\phi}_s,s,X^n)b^2(s,X^n)+\int_Ef(s,X^,u)
\\
&&\times\big(\exp[\Lambda_{T,\phi}(\dot{\phi}_s,s,X^n)f(s,X^n,u)]-1\big)q(du)
\Big]ds.
\end{eqnarray*}
This implies by the definition of
$\Lambda_{T,\phi}(\dot{\phi}_s,s,X^n)$ (see (2.2) and (2.4)) that
\begin{equation}
\langle
Z^{n,\tau},M^n\rangle_t=\int_0^{t\wedge\tau}Z^{n,\tau}_{s-}(\dot{\phi}_s-
a(s,X^n))ds.
\end{equation}
By Theorem 4.5.2 [2], the process $M^{n,\phi}=(M^{n,\phi}_t)_{t\ge
0}$ which is defined by
\[
M^{n,\phi}_t=M^n_t-\int_0^t(Z^{n,\tau}_{s-})^{-1}d\langle
Z^{n,\tau},M^n\rangle _s,
\]
is a local martingale with respect to $(\mathbb{F}^n,Q^{n,\phi})$.
Together with (4.15) this implies that
\begin{equation}
M^{n,\phi}_t=M^n_t-(\phi_{t\wedge\tau}-x)+\int_0^{t\wedge\tau}a(s,X^n)ds.
\end{equation}
Decomposition (4.9) follows from (4.16) and (4.13).

We prove now (4.10). Let us first calculate the quadratic
variation process
$$
[M^{n,\phi},M^{n,\phi}].
$$
Let $M^{n,\phi,c}$ and $M^{n,c}$ be the continuous martingale
components on $M^{n,\phi}$ and $M^n$. The processes $\langle
M^{n,\phi}\rangle$ and $\langle M^n\rangle$ are
$Q^{n,\phi}$-indistinguishable ([2], Theorem 4.5.2), and hence by
(4.4) we have
\begin{equation}
\langle M^{n,\phi,c}\rangle_t=\frac{1}{n}\int_0^tb^2(s,X^n)ds.
\end{equation}
Next, by (4.16) and (4.14)
\begin{equation}
(\triangle M^{n,\phi}_s)^2=n^{-2}\int_Ef^2(s,X^n,u)p^n(\{s\},du).
\end{equation}
Therefore the quadratic variation process
$[M^{n,\phi},M^{n,\phi}]= ([M^{n,\phi},M^{n,\phi}]_t)_{t\ge 0}$ is
\begin{equation}
[M^{n,\phi},M^{n,\phi}]_t=n^{-1}\int_0^tb^2(s,X^0ds+n^{-2}\int_0^t\int_E
f^2(s,X^n,u)p^n(ds,du).
\end{equation}
Now we prove that the increasing process $[M^{n,\phi},M^{n,\phi}]$
is locally integrable with the compensator given by (4.10). Let
$q^{nm\phi}= q^{n,\phi}(ds,du)$ be the compensator of $p^n(ds,du)$
with respect to $(\mathbb{F}^n,Q^{n,\phi})$. By Theorem 4.5.1 [2],
\begin{equation}
q^{n,\phi}(ds,du)=Y(s,u)q^n(ds,du), \ Q^{n,\phi}-\text{a.s.},
\end{equation}
where
\[
Y(s,u)=M^P_{p^n}\Big(1+\frac{\triangle
Z^{n,\phi}}{Z^{n,\phi}}\Big|\mathcal{P} (\mathbb{F}^n\Big)(s,u),
\]
and $M^P_{p^n}(|\mathcal{P}(\mathbb{F}^n)$ is the conditional
expectation of Doleans-Dade measure
$$
M^P_{p^n}=M^P_{p^n}(d\omega,dt,du)=P(\omega)p^n(dt,du)
$$
with respect to $\mathcal{P}(\mathbb{F}^n)$. The definition of
$Z^{n,\tau}$ yields
\[
1+\frac{\triangle
Z^{n,\tau}_s}{Z^{n,\tau}_{s-}}=\exp\Big[I(s\le\tau)
\Lambda_{T,\phi}(\dot{\phi}_s;s,X^n)\int_Ef(s,X^n,u)p^n(\{s\},du)\Big].
\]
It is easy to deduce that $M^P_{p^n}$-a.s.
\begin{equation}
Y(s,u)=\exp[I(s\le\tau)\Lambda_{T,\phi}(\dot{\phi}_s;,s,X^n)f(s,X^n,,u)].
\end{equation}
(4.20) and (4.21) imply that
\begin{equation}
q^{n,\phi}(ds,du)=n\exp[I(s\le\tau)\Lambda_{T,\phi}(\dot{\phi}_s;,s,X^n)
f(s,X^n,,u)]q(du)ds.
\end{equation}
Consider the following increasing process (cf. (4.22))
\begin{eqnarray*}
\alpha_t&=&n^{-2}\int_0^t\int_Ef^2(s,X^n,u)q^{n,\phi}(ds,du)
\\
&=&n^{-2}\int_0^t\int_Ef^2(s,X^n,u)\exp[I(s\le\tau)
\Lambda_{T,\phi}(\dot{\phi}_s;,s,X^n)f(s,X^n,,u)]q(du)ds.
\end{eqnarray*}
In view of {\bf I} and {\bf III(i)}, $\alpha_t<\infty$ $P$-a.s.
for all $t>0$. Therefore, by the equivalence of $Q^{n,\phi}$ and
$P$, $\alpha_t <\infty$ $Q^{n,\phi}$-a.s. The process
$(\alpha_t)_{t\ge 0}$, being increasing and continuous, is locally
integrable (with $\sigma_m=\inf\{t:\alpha_t\ge m\}$, $m\ge 1$ as a
localization sequence). Hence, $(\alpha_t)_{t\ge 0}$ is the
compensator of the process
\[
\Big(n^{-2}\int_0^t\int_Ef^^2(s,X^n,u)p^n(ds,du)\Big)_{t\ge 0}
\]
(with respect to $(\mathbb{F}^n,Q^{n,\phi})$). Then by (4.19)
$[M^{n,\phi},M^{n,\phi}]$ is locally integrable and its
compensator is given by (4.10). This also means that $M^{n,\phi}$
is a local square integrable martingale (see the Burkholder-Gandi
inequality for $p=2$, [2], ch. 1), and so its predictable
quadratic variation process is the compensator of its quadratic
variation process (Theorem 1.8.1 [2]).

\section{\bf Ergodic properties}

Let $T>0$ and let $\phi=(\phi_t)_{t\ge 0}$ be an absolutely
continuous function with $\phi_0=x \ (\equiv X^n_0)$ and bounded
everywhere derivative $\dot{\phi}_t$ on $[0,T]$. We also assume
that $I_T(\phi)<\infty$ and that condition {\bf III(i)} holds.
Define the measure $Q^{n,\phi}$ by (4.7), where $\overline{
\lambda}(t)$ is from (4.6) and $\tau$ from (4.4) with
$\gamma=\delta_{T,\phi}$. $X^n$ is a semimartingale both under $P$
and $Q^{n,\phi}$ with the respective decomposition (4.13) and
(4.9).

In this section we prove that
\begin{equation}
\lim_nP\Big(\sup_{t\le T}|X^n_t-Y_t|>\varepsilon\Big)=0, \ \forall
\varepsilon>0
\end{equation}
and
\begin{equation}
\lim_nQ^{n,\phi}\Big(\sup_{t\le T}|X^n_t-\phi_t|>\delta\Big)=0, \
\forall \delta>0,
\end{equation}
where $(Y_t)_{0\le t\le T}$ is a (unique) solution of the equation
\begin{equation}
\dot{Y}_t=a(t,Y), \ \ Y_0=x.
\end{equation}

\begin{lemma}
Assume that conditions {\bf I} and {\bf II} hold, and {\rm (5.3)}
has a unique solution. Then {\rm (5.1)} holds for all $T>0$.
\end{lemma}

\medskip
{\it Proof.} We denote $Q^{X^n}$ the distribution of $X^n$ under
$P$, in other words $Q^{X^n}$ is the measure on $(D,\mathcal{D})$
defined as
\[
Q^{X^n}(\Gamma)=P(X^n\in\Gamma), \ \ \Gamma\in\mathcal{D}.
\]
We have shown in Section 3 that the sequence $(Q^{X^n})_{n\ge 1}$
is exponentially tight, so since exponential tightness obviously
implies tightness the sequence $(Q^{X^n})$ is tight.

Now we prove that
\begin{equation}
M^{n*}_T\xrightarrow{P}0 \ (n\to\infty).
\end{equation}
By the Lenglart-Rebolledo inequality (see, e.g. [2], ch. 1), (5.4)
will follow from
\begin{equation}
\langle M^n\rangle_T\xrightarrow{P}0 \ (n\to\infty)
\end{equation}
(Problem 1.9.2 [2]). (4.14) implies by condition {\bf I} that
\begin{eqnarray*}
\langle
M^n\rangle_t&=&n^{-1}\int_0^tb^2(s,X^n)ds+n^{-1}\int_0^t\int_E
f^2(s,X^n,u)q(du)ds
\\
&\le&2n^{-1}l^2_T\int_0^t(1+(X^{n*}_s)^2)ds+2n^{-1}l^2_T\int_Eh^2(u)q(du)
\int_0^t(1+(X^{n*}_s)^2)ds.
\end{eqnarray*}
Then there exists $k$ such that
\[
\langle M^n\rangle_T\le kn^{-1}(1+(X^{n*}_T)^2).
\]
Hence
\[
\langle M^n\rangle_TI(X^{n*}_T\le C)\xrightarrow{P}0 \
(n\to\infty), \ \forall C>0,
\]
and (5.5) follows since by the tightness of $Q^{X^n}$
\[
\mathop{\varlimsup}\limits_nP(X^{n*}_T>C)\to 0 \ (C\to\infty).
\]

Denote for $X\in D$
\[
M_t(X)=X_t-X_0-\int_0^ta(s,X)ds.
\]
Then by (5.4)
\begin{equation}
Q^{X^n}\Big(\sup_{t\le
T}|M_t(X)|>\varepsilon\Big)=P(M^{n*}_T>\varepsilon) \to 0 \
(n\to\infty).
\end{equation}
The sequence $Q^{X^n}$ is tight, so it os relatively compact. Let
a subsequence $Q^{X^n}$ converge weakly to a probability measure
$Q'$ on $(D,\mathcal{D}$. Then by (5.6)
\[
\lim_{n'}Q^{X^{n'}}\Big(\sup_{t\le T}|M_t(X)|>\varepsilon\Big)=0.
\]
Since $\sup\limits_{t\le T}|\triangle X^n_s|\xrightarrow{P}0 \
(n\to\infty)$, $Q^{X^n}$ is $C$-tight [3]. Therefore, by condition
{\bf II}, $\sup\limits_ {t\le T}|M_t(X)|=0$ $Q'$-a.s. for every
$T>0$, i.e. $X=(X_t)_{t\ge 0}$ is a solution of (5.3) $Q'$-a.s.
The solution being unique, $Q'$ is the Dirac measure concentrated
on $Y$ for any sequence $(n')$.

\begin{theorem}
Let conditions {\bf I} and {\bf III(i)} hold. Then for all $T>0$
we have {\rm (5.2)}.
\end{theorem}

\medskip
{\it Proof.} Show that for $\delta\le\delta_{T,\phi}$
\begin{equation}
Q^{n,\phi}\Big(\sup_{t\le T}|X^n_t-\phi_t|>\delta\Big)\le
2Q^{n,\phi}\Big(\sup_{t\le T}|X^n_t-\phi_t|\ge\delta\Big).
\end{equation}
Since
\[
\{\tau<T\}\subseteq\Big\{\sup_{t\le\tau}|X^n_t-\phi_t|\ge\delta_{T,\phi}\Big\},
\]
we have
\begin{eqnarray*}
Q^{n,\phi}\Big(\sup_{t\le T}|X^n_t-\phi_t|>\delta\Big)&\le&
Q^{n,\phi}\Big(\sup_{t\le T}|X^n_t-\phi_t|>\delta,\tau=T\Big)
\\
&&+Q^{n,\phi}\Big(\sup_{t\le T}|X^n_t-\phi_t|\ge\delta\Big)
\end{eqnarray*}
and (5.7) follows. So it suffices for us to show that for
$\delta\le\delta_{T,\phi}$
\begin{equation}
\lim_nQ^{n,\phi}\Big(\sup_{t\le T}|X^n_t-\phi_t|\ge\delta\Big)=0.
\end{equation}
Wee have from (4.9) that
\[
\sup_{t\le\tau}|X^n_t-\phi_t|=M^{n,\phi*}_\tau,
\]
and, so all we need to prove is
\begin{equation}
Q^{n,\phi}(M^{n,\phi*}_\tau\ge \delta)=0.
\end{equation}
By the Lenglart-Rebolledo inequality,
\begin{equation}
Q^{n,\phi}(M^{n,\phi*}_\tau\ge \delta)\le
\beta/\delta^2+Q^{n,\phi}(\langle M^{n,\phi}\rangle
_\tau\ge\beta).
\end{equation}
Bt Theorem 4.1 and in view of conditions {\bf I} and {\bf III(i)},
\begin{eqnarray}
\langle
M^{n,\phi}\rangle_t&\le&2n^{-1}l^2_T(1+(X^{n*}_{t-})^2)+2n^{-1}l^2_TT\int_Eh^2(u)
\notag\\
&&\times\exp\{rh(u)(1+X^{n*}_{t-})q(du)(1+(X^{n*}_{t-})^2).
\end{eqnarray}
Since $X^{n*}_{t-}\le \phi^*_t+\delta_{T,\phi}$, by (5.1) there
exists $k$ such that $\langle M^{n,\phi}\rangle_\tau\le k/n,$
which implies in view of (5.10) that
\[
\mathop{\varlimsup}\limits_nQ^{n,\phi}(M^{t,\phi*}_\tau\ge\delta)\le
\beta/\delta^2\to 0 \ (\beta \to 0),
\]
which proves (5.9).

\section{\bf Upper bound}

1. We define $I_T(\phi)$ as in Section 2. In this section we prove
the following

\begin{theorem}
Let condition {\bf I} and {\bf II} hold. Then for all $\phi\in D$
with $\phi_0=x \ (\equiv X^n_0)$
\[
\mathop{\varlimsup}\limits_{\delta\to
0}\mathop{\varlimsup}\limits_nn^{-1}\log P\Big( \sup_{t\le
T}|X^n_t-\phi_t|\le\delta\Big)\le -I_T(\phi).
\]
\end{theorem}

For the proof we will need an auxiliary result.

2. Consider a step function $(\lambda(t))_{t\ge 0}$ of the form
$(0\le t_0<t_1<\ldots<t_k)$
\begin{equation}
\lambda(t)=\sum_{i=1|}^k\lambda_iI_{]t_{i-1},t_i]}(t).
\end{equation}
For $\phi\in D$ set by definition
\begin{equation}
\int_0^t\lambda(s)d\phi_s=\sum_{i=1}^k\lambda_i[\phi_{t\wedge
t_i}-\phi_{t\wedge t_{i-1}}],
\end{equation}
and denote
\begin{equation}
\widetilde{I}_T(\phi)=\sup\Big[\int_0^T\lambda(t)d\phi_t-\int_0^TG(\lambda(t);t,\phi)dt\Big],
\end{equation}
with ``sup'' taken over all $\lambda(t)$ with representation
(6.1).

\begin{lemma}
If $\phi=(\phi_t)_{t\ge 0}\in D$ and $\phi_0=x \ (\equiv X^n_0)$,
then for all $T>0$
\[
I_T(\phi)=\widetilde{I}_T(\phi).
\]
\end{lemma}

\medskip
\noindent {\it Proof.} On first step is that
\begin{equation}
\widetilde{I}_T(\phi)=\infty
\end{equation}
if $\phi$ is not absolutely continuous on $[0,T]$.

By the definition of absolutely continuity [22], for such $\phi$
we can choose $\varepsilon>0$, so that for any $\gamma>0$ there
exists non overlapping intervals $(t'_i,t''_i)$, $i\ge 1$ on
$[0,T]$ satisfying
\begin{equation}
\sum_i(t''_i-t'_i)\le\gamma, \
\sum_i|\phi_{t''_i}-\phi_{t'_i}|>\varepsilon.
\end{equation}
We set an integer $N$
\[
\gamma(N)=[Nl_T(1+\phi^*_T)+(Nl_T)^2(1+(\phi^*_T)^2)+K(Nl_T(1+\phi^*_T))]^{-1},
\]
where $K(\lambda)$ and $l_T$ are defined in condition {\bf I}, and
choose $(t'_i,t''_i)$ for each $N$ so that (6.5) holds with
$\gamma=\gamma(N)$. Let
\begin{equation}
\lambda^0(t)=N\sum_{i=1}^k\text{sign} \
(\phi_{t''_i}-\phi_{t'_i})I_{]t'_i,t''_i]} (t).
\end{equation}
By condition {\bf I}, $|G(\lambda^0(t);t\phi)|\le \gamma^{-1}(N)$
for all $t\le T$, and hence
\begin{equation}
\int_0^T|G(\lambda^0(t);t,\phi)|dt=\int_0^TI(\lambda^0(t);t,\phi)|dt\le
\gamma^{-1}(N)\sum_{i=1}^k(t''_i-t'_i)\le 1.
\end{equation}
On the other hand,
\[
\int_0^T\lambda^0(t)d\phi_t=N\sum{i=1}^k|\phi_{t''_i}-\phi_{t'_i}|>N\varepsilon.
\]
It then follows by the definition of $\widetilde{I}_T(\phi)$ and
(6.2), that
\begin{eqnarray*}
\widehat{I}_T(\phi)&\ge&\int_0^T\lambda^0(t)d\phi_t-\int_0^TG(\lambda^0(t);t,
\phi)dt
\\
&\ge&\int_0^T\lambda^0(t)d\phi_t-\int_0^T|G(\lambda^0(t);t,
\phi)|dt\ge N\varepsilon-1\to\infty, \ N\to\infty
\end{eqnarray*}
which yields (6.4).

Now let $\phi$ be absolutely continuous on $[0,T]$. We show that
\begin{equation}
\widehat{I}_T(\phi)=\int_0^TH(\dot{\phi}_t;t,\phi)dt.
\end{equation}
First, we have by the definition of $H(\dot{\phi}_t;t,\phi)$ (see
(2.3))
\[
\widehat{I}_T(\phi)\le\int_0^TH(\dot{\phi}_t;t,\phi)dt.
\]
So, it is left tto prove
\begin{equation}
\widehat{I}_T(\phi)\ge\int_0^TH(\dot{\phi}_t;t,\phi)dt.
\end{equation}
It is rather obvious that (6.9) holds, provided we can choose for
each $c>0$ a finite function $\lambda^c(t)$ satisfying
\begin{equation}
\lambda^c(t)\dot{\phi}_t-G(\lambda^c(t);t,\phi)\ge [(c\wedge
H(\dot{\phi}_t;t, \phi))-1/c]\vee 0
\end{equation}
for almost all $t\le T$ taking
$\tilde{\lambda}^c(t)=\lambda^c(t)I( |\lambda^c(t)|\le c))$.
Indeed there exists a sequence $\lambda_n(t)$, $n\ge 1$ of
step-functions of the form (6.1) for which
\begin{equation}
|\lambda_n(t)|\le \sup_{t\le T}|\tilde{\lambda}^c(t)|, \ t\le T,
\end{equation}
\begin{equation}
\lim_n\int_0^T|\lambda_n(t)-\tilde{\lambda}^c(t)|dt=0
\end{equation}
(Lemma 4.4, ch. 4 in [23]). By condition {\bf I} and (6.11)
\begin{equation}
\big|[\tilde{\lambda}^c(t)\dot{\phi}_t-G(\tilde{\lambda}^c(t);t,\phi))]-
[\lambda_n(t)\dot{\phi}_t-G(\lambda_n(t);t,\phi))]\big|
\le|\lambda_n(t)-\tilde{\lambda}^c(t)|(|\dot{\phi}_t|+r),
\end{equation}
where $r$ depends on $\phi$, on the parameters in condition {\bf
I} and on $c$. So, provided (6.10) holds. we have using the
definition $\widetilde{I}_T(\phi)$
\begin{eqnarray}
\widetilde{I}_T(\phi)&\ge&
\int_0^T\lambda_n(t)d\phi_t-\int_0^TG(\lambda_n(t); t,\phi)dt
\notag\\
&\ge&\int_0^T\tilde{\lambda}^c(t)d\phi_t-\int_0^TG(tilde{\lambda}^c(t);
t,\phi)dt
\\
&&-\int_0^T|\lambda_n(t)-\tilde{\lambda}^c(t)|(|\dot{\phi}_t|+r)dt
\notag\\
&\ge&\int_0^TI\{|\lambda^c(t)|\le c\}\{[(c\wedge
H(\dot{\phi}_t;t\phi))-1/c] \vee 0\}dt
\notag\\
&&-\int_0^T|\lambda_n(t)-\tilde{\lambda}^c(t)|(|\dot{\phi}_t|+r)dt.
\end{eqnarray}
We next show that
\begin{equation}
\lim_n\int_0^T|\lambda_n(t)-\tilde{\lambda}^c(t)|(|\dot{\phi}_t|+r)dt=0.
\end{equation}
(6.15) holds by the inequality (see (6.11))
\begin{eqnarray*}
\int_0^T|\lambda_n(t)-\tilde{\lambda}^c(t)|(|\dot{\phi}_t|+r)dt&\le&(N+r)
\int_0^T|\lambda_n(t)-\tilde{\lambda}^c(t)|dt
\\
&&+2\sup_{t\le
T}|\tilde{\lambda}^c(t)|\int_0^TI(|\dot{\phi}_t|>N)dt,
\end{eqnarray*}
since in view of (6.12) and by the absolute continuity of $\phi$
``$\mathop{\varlimsup}\limits_N\mathop{\varlimsup}\limits_n$'' of
the latter expression is 0.

Thus, (6.14) and (6.15) imply that under (6.10) for any $c>0$
\[
\widetilde{I}_T(\phi)\ge\int_0^TI\{|\tilde{\lambda}^c|\le
c\}\{[(c\wedge H(\dot{\phi}_t;t,\phi))-1/c]\vee 0\}dt.
\]
Since $H(\dot{\phi}_t;t,\phi)$ is nonnegative, we have
\[
\lim_{c\to\infty}\int_0^TI\{|\tilde{\lambda}^c|\le c\}\{[(c\wedge
H(\dot{\phi}_t;t,\phi))-1/c]\vee
0\}dt=\int_0^tH(\dot{\phi}_t;,t\phi)dt
\]
which proves (6.9).

It is left to find $\tilde{\lambda}^c(t)$ meeting (6.10). We
denote
\[
G_0(\lambda;t,\phi)=G(\lambda;t,\phi)-\lambda a(t,\phi)
\]
and
\begin{eqnarray*}
&&
U'(t,\lambda)=\lambda|\dot{\phi}_t-a(t,\phi)|-G_0(\lambda;t,\phi),
\\
&&
U''(t,\lambda)=-\lambda|\dot{\phi}_t-a(t,\phi)|-G_0(\lambda;t,\phi).
\end{eqnarray*}
Introduce the functions
\begin{eqnarray*}
&& \lambda'(t)=\inf\{\lambda\ge 0:U'(t,\lambda)\ge [(c\wedge
H(\dot{\phi}_t;t,\phi) )-1/c]\vee 0\}
\\
&& \lambda''(t)=\inf\{\lambda\ge 0:U''(t,-\lambda)\ge [(c\wedge
H(\dot{\phi}_t;t,\phi) )-1/c]\vee 0\},
\end{eqnarray*}
($\inf(\varnothing)=\infty$). Then
\[
\tilde{\lambda}^c(t)=\lambda'(t)I(\dot{\phi}_t>a(t,\phi))\lambda''(t)I(
\dot{\phi}_t<a(t,\phi)).
\]

\medskip
{\it Proof of Theorem 6.1}. Let $\lambda(t)$ be of the form (6.1).
Define a positive local martingale $Z^n=(Z^n_t)_{t\ge 1}$ by
(4.3). By Lemma 4.1, $EZ^n_T=1$ and hence
\begin{equation}
1\ge E\Big[I\Big(\sup_{t\le
T}|X^n_t-\phi_t|\le\delta\Big)Z^n_T\Big].
\end{equation}
The expression under the  expectation is less than
\begin{multline*}
I\Big(\sup_{t\le
T}|X^n_t-\phi_t|\le\delta\Big)\exp\Big\{n\Big[\int_0^T
\lambda(t)\dot{\phi}_tdt-\int_0^TG(\lambda(t);t,\phi)dt\Big]
\\
\qquad\qquad-n\Big[\Big|\int_0^T\lambda(t)d(X^n_t-\phi_t)\Big| -
\Big|\int_0^TG(\lambda(t);t,\phi)dt-\int_0^TG(\lambda(t);t,X^)dt\Big|\Big]\Big\}.
\end{multline*}
Since $\lambda(t)$ is piece wise constant, we have on
$\{\sup_{t\le T}|X^n_t- \phi_t|\le\delta\}$
\[
\Big|\int_0^T\lambda(t)d(X^n_t-\phi_t)\Big|\le 2\sup_{t\le
}|\lambda(t)|\delta.
\]
By conditions {\bf I} and {\bf II}, $\int_0^TG(\lambda(t);t,X)dt$
is continuous in $X\in D$ at each $\phi\in C$, i.e. if $\sup_{t\le
T}|X^{(k)}_t-\phi_t|\to 0 \ (k\to\infty)$ for $X^{(k)}\in D$,
$k\ge 1$, then
\[
\int_0^TG(\lambda(t);t,X^{k)})dt\to\int_0^TG(\lambda(t);t,\phi)dt
\ (k\to\infty).
\]
It follows that for each for each $\varepsilon>0$ there exists
$\delta( \varepsilon,\phi)$ such that for
$\delta<\delta(\varepsilon,\phi)$ and $\sup_{t\le
T}|X^n_t-\phi_t|\le\delta$
\[
\Big|\int_0^TG(\lambda(t);t,X^{k)})dt-\int_0^TG(\lambda(t);t,\phi)dt\Big|\le
\varepsilon.
\]
Hence for $\delta<\delta(\varepsilon,\phi)$
\begin{eqnarray*}
&& I\Big(\sup_{t\le T}|X^n_t-\phi_t|\le\delta\Big)Z^n_T
\\
&\ge& I\Big(\sup_{t\le T}|X^n_t-\phi_t|\le\delta\Big)
\\
&&\times \exp\Big\{n\Big[\int_0^T
\lambda(t)\dot{\phi}_tdt-\int_0^TG(\lambda(t);t,\phi)dt-2\sup_{t\le
T} |\lambda(t)|\delta-\varepsilon\Big]\Big\}.
\end{eqnarray*}
This and (6.16) imply since $\lambda(t)$ and $\varepsilon$ are
arbitrary, that
\begin{eqnarray*}
&& \mathop{\varlimsup}\limits_{\delta\to
0}\mathop{\varlimsup}\limits_nn^{-1}\log P\Big(\sup_{t\le
T}|X^n_t-\phi_t|\le\delta\Big)
\\
&\le&-\sup\Big[\int_0^T\lambda(t)\dot{\phi}_tdt-\int_0^TG(\lambda(t);t,\phi)dt
\Big],
\end{eqnarray*}
with ``sup'' taken over all step-functions $\lambda(t)$ of the
form (6.1).

The theorem follows by Lemma 6.1.

\section{\bf Lower bound}

\begin{lemma}
Let $\phi=(\phi_t)_{t\ge 0}$ be absolutely continuous on $[0,T]$
with $\phi_0=x \ (\equiv X^n_0)$ and bounded almost everywhere
derivative: $|\dot{\phi}_t|\le N$. If conditions {\bf I}, {\bf II}
and {\bf III} hold then
\[
\mathop{\varliminf}\limits_{\delta\to
0}\mathop{\varliminf}\limits_nn^{-1}\log P\Big(\sup_{t\le
T}|X^n_t-\phi_t|\le\delta\Big)\ge -I_T(\phi).
\]
\end{lemma}

\medskip
{\it Proof.} Since $\dot{\phi}_t$ is bounded, so by condition {\bf
III}
\begin{equation}
I_T(\phi)<\infty.
\end{equation}
Define $Q^{n,\phi}$ by (4.7) with $\lambda(t)$ from (4.6) and
$\tau$ from (4.4) with $\gamma=\delta_{T,\phi}$. Then by (4.8)
\begin{equation}
P\Big(\sup_{t\le T}|X^n_t-\phi_t|\le\delta\Big)=\int_\Omega I\Big(
\sup_{t\le
T}|X^n_t-\phi_t|\le\delta\Big)(Z^{n,\tau}_\infty)^{-1}dQ^{n,\phi},
\end{equation}
where
\begin{equation}
Z^{n,\tau}_\infty=\exp\Big\{n\Big[\int_0^\tau
\Lambda_{T,\phi}(\dot{\phi}_t;t,X^n) dX^n_t-\int_0^\tau
G(\Lambda_{T,\phi}(\dot{\phi}_t;t,X^n);t,X^n)dt\Big]\Big\}.
\end{equation}
We estimate the right hand side of (7.2). Note that by condition
{\bf III}
\begin{eqnarray}
&& \int_0^\tau\{\Lambda_{T,\phi}(\dot{\phi}_t;t,X^n)\dot{\phi}_t-
G(\Lambda_{T,\phi}(\dot{\phi}_t;t,X^n);t,X^n)\}dt
\notag\\
&=& \int_0^\tau H(\dot{\phi}_t;t\phi)dt\le \int_0^T
H(\dot{\phi}_t;t\phi)dt =I_T(\phi).
\end{eqnarray}
Let $M^{n,\phi}=(M^{n,\phi}_t)_{t\ge 0}$ be local square
integrable martingale (under $Q^{n,\phi}$) from Theorem 4.1.

By condition {\bf I} and Theorem 5.1, and since $|\dot{\phi}_t|$
is bounded, $\Lambda_{T,\phi}(\dot{\phi}_t;t,X^n)$ is uniformly
bounded for almost all $t\le\tau$ and so the integral
\[
\int_0^\tau \Lambda_{T,\phi}(\dot{\phi}_t;t,X^n)dM^{n,\phi}_t
\]
ia well defined. Then (7.3) and (7.4) yield
\begin{equation}\label{7.5}
\begin{aligned}
(Z^{n,\tau}_\infty)^{-1}&\ge\exp\Big\{-n\Big[I_T(\phi)+
\int_0^\tau \Lambda_{T,\phi}(\dot{\phi}_t;t,X^n)dM^{n,\phi}_t
\\
&+\int_0^T|\dot{\phi}_t||\Lambda_{T,\phi}(\dot{\phi}_t;t,X^n)-
\Lambda_{T,\phi}(\dot{\phi}_t;t,\phi)|dt
\\
&+\int_0^T|G(\Lambda_{T,\phi}(\dot{\phi}_t;t,X^n);t,X^n)-
G(\Lambda_{T,\phi}(\dot{\phi}_t;t,\phi);t,\phi|dt\Big]\Big\}.
\end{aligned}
\end{equation}
Conditions {\bf I} - {\bf III} and boundedness of $\dot{\phi}_t$
imply that for every $\varepsilon>0$ there exists
$\delta(\varepsilon,T,\phi)\le \delta_{T, \phi}$ (which depends on
$\varepsilon,T$ and $\phi$) such that, provided
\[
\sup_{t\le T}|X^n_t-\phi_t|\le \delta(\varepsilon,T,\phi),
\]
we have
\begin{eqnarray}
&& \int_0^T|\dot{\phi}_t||\Lambda_{T,\phi}(\dot{\phi}_t;t,X^n)-
\Lambda_{T,\phi}(\dot{\phi}_t;t,\phi)|dt
\notag\\
&&+\int_0^T|G(\Lambda_{T,\phi}(\dot{\phi}_t;t,X^n);t,X^n)-
G(\Lambda_{T,\phi}(\dot{\phi}_t;t,\phi);t,\phi|dt\Big]\Big\}
\notag\\
&\le&\varepsilon.
\end{eqnarray}
(7.2) and (7.3), and (7.5) and (7.6) yield for $\delta\le
\delta_{\varepsilon,T, \phi}$:
\begin{eqnarray}
P\Big(\sup_{t\le T}|X^n_t-\phi_t|\le
\delta\Big)&\ge&\exp\big(-n\big[I_T(\phi)+
\varepsilon]\big)\int_\Omega\Big[I\Big(\sup_{t\le
T}|X^n_t-\phi_t|\le \delta\Big)
\notag\\
&&\times\exp\Big(-n\Big|\int_0^\tau
\Lambda_{T,\phi}(\dot{\phi}_t;t,X^n)
dM^{n,\phi}_t\Big|\Big)\Big]dQ^{n,\phi}.
\end{eqnarray}
We have for the integral of the right hand side of (7.7) for
$\beta>0$
\begin{eqnarray}
&& \int_\Omega\Big[I\Big(\sup_{t\le T}|X^n_t-\phi_t|\le
\delta\Big) \exp\Big(-n\Big|\int_0^\tau
\Lambda_{T,\phi}(\dot{\phi}_t;t,X^n)
dM^{n,\phi}_t\Big|\Big)\Big]dQ^{n,\phi}
\notag\\
&\ge&\exp\{-n\beta\}Q^{n,\phi}\Big\{\Big|\int_0^\tau
\Lambda_{T,\phi} (\dot{\phi}_t;t,X^n)dM^{n,\phi}_t\Big|\le\beta,
\sup_{t\le T}|X^n_t-\phi_t|\le \delta\Big\}.
\end{eqnarray}

(7.7) and (7.8) imply that for $\delta\le
\delta(\varepsilon,T,\phi)$
\begin{eqnarray}
&& n^{-1}\log P\Big(\sup_{t\le T}|X^n_t-\phi_t|\le\delta\Big)\ge
-(I_T(\phi)+ \varepsilon+\beta)
\notag\\
&& +n^{-1}\log Q^{n,\phi}\Big\{\Big|\int_0^\tau \Lambda_{T,\phi}
(\dot{\phi}_t;t,X^n)dM^{n,\phi}_t\Big|\le\beta,\sup_{t\le
T}|X^n_t-\phi_t|\le \delta\Big\}.
\end{eqnarray}
By (7.9) and Theorem 5.1, it suffices for us to prove that for any
$\beta>0$
\begin{equation}
Q^{n,\phi}\Big\{\Big|\int_0^\tau \Lambda_{T,\phi}
(\dot{\phi}_t;t,X^n)dM^{n,\phi}_t\Big|> \beta\Big\}=0.
\end{equation}
Since by condition {\bf III}
$|\Lambda_{T,\phi}(\dot{\phi}_t;t,X^n)|\le r$ for almost all
$s\le\tau$, the Lenglart-Rebolledo inequality gives
\begin{eqnarray}
&& Q^{n,\phi}\Big\{\Big|\int_0^\tau \Lambda_{T,\phi}
(\dot{\phi}_t;t,X^n)dM^{n,\phi}_t\Big|> \beta\Big\}
\notag\\
&\le&\alpha/\beta^2+Q^{n,\phi}\Big\{\int_0^\tau \Lambda^2_{T,\phi}
(\dot{\phi}_t;t,X^n)d\langle M^{n,\phi}\rangle_t\ge \beta\Big\}
\notag\\
&\le& \alpha/\beta^2+Q^{n,\phi}\{\langle M^{n,\phi}\rangle_\tau\ge
\alpha/r^2\}.
\end{eqnarray}
While proving Theorem 5.1 we saw that $\langle
M^{n,\phi}\rangle_\tau\le k/n$, where $k$ does not depend on $n$.
Taking on the right side of (7.11)
``$\lim_{\alpha\to0}\varlimsup_n$'' we obtain (7.10).

\medskip
\begin{theorem}
Let $\phi=(\phi_t)_{t\ge 0}\in D$ and $\phi_0=x \ \equiv X^n_0)$.
Under {\bf I} - {\bf IV} we have
\begin{equation}
\lim_{\delta\to 0}\mathop{\varliminf}\limits_nn^{-1}\log
P\Big(\sup_{t\le T} |X^n_t-\phi_t|\le\delta\Big)\ge -I_T(\phi).
\end{equation}
\end{theorem}

\medskip
{\it Proof.} Obviously we can assume (7.1), in particular $\phi$
is absolutely continuous on $[0,T]$. Define
$\phi^n=(\phi^N_t)_{t\ge 0}$ as in condition {\bf IV}:
\begin{equation}
\phi^N_t=x+\int_0^tI(|\dot{\phi}_s|\le N)\dot{\phi}_sds.
\end{equation}
Applying Lemma (7.1) to $\phi^N$ we have for all $T>0$
\begin{equation}
\lim_{\delta\to 0}\mathop{\varliminf}\limits_nn^{-1}\log
P\Big(\sup_{t\le T} |X^n_t-\phi^N_t|\le\delta\Big)\ge
-I_T(\phi^N).
\end{equation}
Since
\[
\lim_N\sup_{t\le T}|\phi_t-\phi^N_t|=0,
\]
\[
P\Big(\sup_{t\le T}|X^n_t-\phi_t|+\sup_{t\le T}|\phi_t-\phi^N_t|
\le\delta\Big)\ge P\Big(\sup_{t\le T}|X^n_t-\phi_t|\le\delta\Big).
\]
we have from (7.14) that
\[
\lim_{\delta\to 0}\mathop{\varliminf}\limits_nn^{-1}\log
P\Big(\sup_{t\le T}|X^n_t-\phi_t|\le\delta\Big)\ge
-\mathop{\varliminf}\limits_N I_T(\phi^N),
\]
which implies that the assertion will follow if
\begin{equation}
\mathop{\varliminf}\limits_NI_T(\varphi^N)\le I_T(\varphi).
\end{equation}
using the definition of $I_T(\phi)$ (see Section 2) we have
\begin{eqnarray}\
I_T(\phi^N)&=&\int_0^TH(\dot{\phi}^N;t,\phi^N)dt
\notag\\
&=&\int_0^TH(\dot{\phi}^N;t,\phi^N)I(|\dot{\phi}_t|\le N)dt+
\int_0^TH(\dot{\phi}^N;t,\phi^N)I(|\dot{\phi}_t|> N)dt
\notag\\
&=&\int_0^TH(\dot{\phi};t,\phi^N)I(|\dot{\phi}_t|\le N)dt+
\int_0^TH(0;t,\phi^N)I(|\dot{\phi}_t|> N)dt
\notag\\
&\le&\int_0^TH(\dot{\phi};t,\phi^N)dt+
\int_0^TH(0;t,\phi)I(|\dot{\phi}_t|> N)dt
\notag\\
&&+\int_0^T|H(0;t,\phi)-H(0;t,\phi^N)|dt.
\end{eqnarray}
Taking into account conditions {\bf I}, {\bf II} and {\bf III} it
is easy to to see that the last two terms on the right hand side
of (7.16) go to 0 as $N\to\infty$. Then (7.15) holds if
\[
\mathop{\varliminf}\limits_N\int_0^TH(\dot{\varphi}_t;t,\phi^N)dt\le
I_T (\varphi).
\]
By conditions {\bf I}, {\bf II} and {\bf III} we have as
$N\to\infty$ for any $m>0$
\begin{eqnarray*}
\int_0^TH(\dot{\varphi}_t;t,\phi^N)I(|\dot{\phi}_t|\le m)dt&\to&
\int_0^TH(\dot{\varphi}_t;t,\phi)I(|\dot{\phi}_t|\le m)dt
\notag\\
&\le&\int_0^TH(\dot{\varphi}_t;t,\phi)dt=I_T(\phi).
\end{eqnarray*}
Hence, it suffices that
\[
\lim_m\mathop{\varlimsup}\limits_N
\int_0^TH(\dot{\varphi}_t;t,\phi^N)I(|\dot{\phi}_t|> m)dt=0,
\]
which is a consequence of condition {\bf IV} and the Lebesgue
dominated convergence theorem.

\medskip
\begin{theorem}
Let the conditions of Theorem 2.2 hold and $\varphi$ be the
solution of (2.9). Then for all $T>0$
\[
\lim_{\delta\to 0}\mathop{\varliminf}\limits_nn^{-1}\log
P\Big(\sup_{t\le T} |X^n_t-\phi_t|\le\delta\Big)\ge
-I_T(\phi)=q(E)T.
\]
\end{theorem}

\medskip
{\it Proof.} We first show that
\begin{equation}
I_T(\phi)=q(E)T.
\end{equation}
Since $\phi$ solves (2.9), we have (see (2.1))
\[
\lambda\dot{\phi}_t-G(\lambda;t,\phi)=\int_E(1-e^{\lambda
f(t,\phi,u})q(du)
\]
and so by (2.3) $H(\dot{\phi}_t;t,\phi)=q(E)$, which gives (7.17).

We now prove the required. For $c>0$ consider the equation
\[
\dot{\phi}^c_t=a(t,\phi^c)+\int_E(e^{-cf(t,\phi,u)}-1)f(t,\phi^c,u)q(du).
\]
It is not difficult to see using conditions {\bf I} and {\bf II}
that the solution $\phi^c=(\phi^c_t)_{t\ge 0}$ exists for $t\in
[0,\infty)$ and it suffices the conditions of Lemma 7.1. Then
\begin{equation}
\lim_{\delta\to 0}\mathop{\varliminf}\limits_nn^{-1}\log
P\Big(\sup_{t\le T} |X^n_t-\phi^c_t|\le\delta\Big)\ge
-I_T(\phi^c),
\end{equation}
where
\[
I_T(\phi^c)=\int_0^T\int_E[1-(c+1)e^{-cf(t,\phi,u)}]q(du)dt.
\]
Obviously
\begin{equation}
I_T(\phi^c)\to q(E)T \ (c\to\infty).
\end{equation}
The inequality
\[
P\Big(\sup_{t\le T}|X^n_t-\phi_t|+\sup_{t\le T}|\phi_t-\phi^N_t|
\le\delta\Big)\ge P\Big(\sup_{t\le T}|X^n_t-\phi_t|\le\delta\Big),
\]
and (7.18) and (7.19) imply that it suffices to prove that
\begin{equation}
\sup_{t\le T}|\phi_t-\phi^c_t|\to 0 \ (c\to\infty), \ \forall \
T>0.
\end{equation}

It is easy to see applying the Arzela-Ascoli theorem that the
family $(\varphi^c_t)_{0\ge t\le T}$, $c>0$ is relatively compact
in $C_{[0,T]}$. Since by conditions {\bf I} and {\bf II} any
subsequential limit of $\phi^c$ as $c\to\infty$ solves (2.9), and
the solution of (2.9) is unique, (7.29) is proved.

\section{\bf Proof of main result}

{\it Proof of Theorem 2.1.} In Theorem 3.2 we proved that the
sequence $X^n,n\ge 1$ is $C$-exponentially tight. Besides, by
Theorem 6.1 for any $\phi\in D$ with $\phi_0=x \ (\equiv X^n_0)$
\[
\mathop{\varlimsup}\limits_{\delta\to
0}\mathop{\varlimsup}\limits_nn^{-1}\log P\Big(\sup_{t\le
T}|X^n_t-\phi_t|\le\delta\Big)\le -I_T(\phi),
\]
where $I_T(\phi)$ is defined in defined in Section 2. On the other
hand, by Theorem 7.1 we have for any $\phi\in D$ with $\phi_0=x \
(\equiv X^n_0)$
\[
\lim_{\delta\to 0}\mathop{\varlimsup}\limits_nn^{-1}\log
P\Big(\sup_{t\le T}|X^n_t-\phi_t|\le\delta\Big)\ge -I_T(\phi),
\]
Therefore for any $\phi\in C$ with $\phi\in D$ with $\phi_0=x \
(\equiv X^n_0)$
\begin{eqnarray}
I_T(\phi)&=&-\mathop{\varlimsup}\limits_{\delta\to
0}\mathop{\varlimsup} \limits_nn^{-1}\log P\Big(\sup_{t\le
T}|X^n_t-\phi_t|\le\delta\Big)
\notag\\
&=&-\mathop{\varliminf}\limits_{\delta\to 0}\mathop{\varliminf}
\limits_nn^{-1}\log P\Big(\sup_{t\le
T}|X^n_t-\phi_t|\le\delta\Big),
\end{eqnarray}
i.e. (1.11) holds with $J_T(\phi)=I_T(\phi)$ for $\phi\in C$ with
$\phi_0=x \ (\equiv X^n_0)$.

Obviously we have (8.1) with $\phi_0\neq x \ (\equiv X^n_0)$ as
well, since in that case each term in (8.1) is  equal to $\infty$.

So according to the scheme in Section 1property $\beta$) holds and
the rate function for $X^n,n\ge 1$ is given in $\gamma$).

The proof of Theorem 2.2 is analogous  with the use of Theorem 7.2
in addition to Theorem 7.1.

\section{\bf Explicit conditions on coefficients}

In this section we give simple sufficient conditions on $a(t,X)$,
$b(t,X)$, $f(t,X,u)$ which imply conditions {\bf III} and {\bf
IV}.

We begin with some definitions. Say that $f^+(t,X,u) \ (=f\vee 0)$
is nondegenerate with respect to $q(du)$ uniformly on $[0,T]$ in a
neighborhood of $\phi=(\phi_t)_{t\ge 0}\in C$, if  there exists
$\delta>0$ and $\gamma>0$ such that for all $X=(X_t)_{t\ge 0}\in
D$ with $\sup_{t\le T}|X_t-\phi_t|\le \delta$
\[
\sup_{t\le T}\int_Ef^+(t,X,u)I(f^+(t,X,u)>\gamma)q(du)>0.
\]

Say that $f^+(t,X,u) \ (=f\vee 0)$ is degenerate with respect to
$q(du)$ uniformly on $[0,T]$ in a neighborhood of
$\phi=(\phi_t)_{t\ge 0}\in C$, if  there exists $\delta>0$ such
that for all $X=(X_t)_{t\ge 0}\in D$ with $\sup_{t\le
T}|X_t-\phi_t|\le \delta$
\[
\sup_{t\le T}\int_Ef^+(t,X,u)q(du)=0.
\]

Similarly the nondegeneracy and degeneracy of $f^-(t,X,u) \
(=-f\wedge 0)$ are defined.

Next, we call $b(t,X)$ nondegenerate uniformly on $[0,T]$ in a
neighborhood of $\phi=(\phi_t)_{t\ge 0}\in C$, if  there exists
$\delta>0$ such that for all $X=(X_t)_{t\ge 0}\in D$ with
$\sup_{t\le T}|X_t-\phi_t|\le \delta$
\[
\inf_{t\le T}b^2(t,X)>0.
\]

\begin{theorem}
Assume that
$\varlimsup_{\lambda\to\infty}K(\lambda)/\lambda<\infty$ {\rm
(}$K(\lambda)$ is defined in condition {\bf I}{\rm )}. Let for ant
$T>0$ and any $\phi=( \phi_t)_{t\ge 0}$ with $I_T(\phi)<\infty$ at
least one of the following conditions hold:

1{\rm }) both $f^+(t,X,u)$ and $f^-(t,X,u)$ are nondegenerate with
respect to $q(du)$ uniformly on $[0,T]$ in a neighborhood of
$\phi$;

2{\rm )} $b(t,X)$ is nondegenerate uniformly on $[0,T]$ in a
neighborhood of $\phi$ and each of $f^+(t,X,u)$ and $f^-(t,X,u)$
is either degenerate or nondegenerate with respect to $q(du)$
uniformly on $[0,T]$ in a neighborhood of $\phi$.

Then conditions {\bf III} and {\bf IV} hold.
\end{theorem}

\smallskip
First wee prove two lemmas.

\begin{lemma}
If for any $T>0$ and any $\phi=(\phi_t)_{t\ge 0}$ with
$I_T(\phi)<\infty$ there exists $\delta>0$ such that for all
$X=(X_t)_{t\ge 0}\in D$ with $\sup_{t\le T}|X_t-\phi_t|\le \delta$
\begin{eqnarray*}
&& c_+=\inf_{t\le T}(b^2(t,X)+\int_E(f^+(t,X,u))^2q(du)>0
\\
&& c_-=\inf_{t\le T}(b^2(t,X)+\int_E(f^-(t,X,u))^2q(du)>0,
\end{eqnarray*}
then condition {\bf III} holds with $\delta_{T,\phi}\delta$.
\end{lemma}

\medskip
{\it Proof.} We have to prove that the equation
\begin{equation}
y=g(\lambda;t,X)
\end{equation}
where (see Section 2)
\begin{equation}
g(\lambda;t,X)=a(t,X)+\lambda b^2(t,X)+\int_Ef(t,X,u)(e^{\lambda
f(t,X,u)}-1) q(du),
\end{equation}
has a solution for all $y\in R$, for all $t\le T$ and $X\in D$
with $\sup_{t\le T}|X_t-\phi_t|\le \delta_{T,\phi}$ for some
$\delta_{T,\phi}>0$ and this solution satisfies {\bf III}.

We begin with proving the existence. Assume that $y\ge a(t,X)$
where $X$ is as in the conditions of the lemma. As for $\lambda>0$
\[
\lambda b^2(t,X)+\int_E f(t,X,u)(e^{\lambda f(t,X,u)}-1)q(du)\ge
\lambda b^2(t,X)+\int_E (f^+(t,X,u))^2q(du),
\]
so $g(\tilde{\lambda};t,X)>y$ for all $t\le T$ where
\[
\tilde{\lambda}=\frac{y-a(t,X)}{h^2(t,X)+\int_E(f^+(t,X,u))^2q(du)}.
\]
Since $g(0;t,X)\le y$ and $g(\lambda;t,X)$ is continuous in
$\lambda$m a solution $\widetilde{\Lambda}_{T,\phi}(y;t,X)$ of
(9.1) exists and satisfies the inequality
\[
0\le\widetilde{\Lambda}_{T,\phi}(y;t,X)\le
\frac{y-a(t,X)}{h^2(t,X)+\int_E(f^+(t,X,u))^2q(du)}.
\]
This solution is unique since under assumptions of the lemma
\begin{equation}
g'_\lambda(\lambda;t,X)=b^2(t,X)+\int_Ef^2(t,X,u)e^{\lambda
f(t,X,u)}q(du))>0.
\end{equation}

The case when $y<a(t,X)$ is considered similarly.

Thus the solution of (9.1) (with respect $\lambda$) is unique and
for all $t\le T$
\begin{equation}
\widetilde{\Lambda}_{T,\phi}(y;t,X)\le\frac{|y-a(t,X)|}{c_+\wedge
c^-}
\end{equation}
for all $X=(X_t)_{t\ge 0}\in D$ with $\sup_{t\le
T}|X_t-\phi_t|\le\delta$. The latter and condition {\bf I} (for
$a(t,X)$) imply $\widetilde{\Lambda}_{T,\phi}(y;t,X)$ meets part
(i) of condition {\bf III} (with $\delta=\delta_{T,\phi}$).

Set
$\Lambda_{T,\phi}(y;t,X)=\widetilde{\Lambda}_{T,\phi}(y;t,X)I(t\le
T) I(\sup_{s\le t}|X_s-\phi_s|\le\delta$.
$\Lambda_{T,\phi}(y;t,X)$ is
$\mathcal{B}(R)\otimes\mathcal{P}(D)$-measurable. Indeed for
$\lambda>0$
\begin{eqnarray*}
&& \{(y,t,X):\Lambda_{T,\phi}(y;t,X)>\lambda\}
\\
&=&\Big\{(y,t,X):g(\lambda;t,X)<y, \ t\le T, \ \sup_{s\le
t}|X_s-\phi_s|\le \delta\Big\}\in
\mathcal{B}(R)\otimes\mathcal{P}(D).
\end{eqnarray*}
The case $\lambda<0$ does not differ.

Finally we prove that $\Lambda_{T,\phi}(y;t,X)$ is $C_{[0,T]}$ in
$X$ at $X=\phi$ for $y\in R$ and all $t\le T$. Let $\sup_{t\le
T}|X^{(k)}_t-\phi_t| \to\infty$ where $X^{(k)}=(X^{(k)}_t)_{t\ge
0}\in D$. Then we can assume that $\sup_{t\le
T}|X^{(k)}_t-\phi_t|\le\delta$ for all $k$ and hence by (9.4) the
sequence $\Lambda^{(k)}=\Lambda_{T,\phi}(y;t,X^{(k)})$ is bounded
in $k$ for all $t\in [0,T]$. Fix $t$ and let $\Lambda^{(k')}$ be a
converging subsequence. Then denoting
\[
\Lambda^{0}=\lim_{k'\to\infty}\Lambda^{(k')}
\]
we have by conditions {\bf I} and {\bf II}, and by (9.2)
\[
y=\lim_{k'\to\infty}g(\Lambda^{(k')};t,X^{(k')})=g(\Lambda^0;t,\phi).
\]
Since the solution of (9.1) for $X=\phi$ is unique, we obtain that
$\Lambda^0=\Lambda_{T,\phi}(y;t,\phi)$.
\begin{lemma}
Let condition {\bf III} hold. Assume that foe any absolutely
continuous function $\phi=(\phi_t)_{t\ge 0}$ with
$I_T(\phi)<\infty$ there exists a nonnegative function
$g(_{T,\phi}(y)$ which satisfies the linear growth condition
\[
g(_{T,\phi}(y)\le C(1+|y|), \ C>0,
\]
and is such that for all $t\le T$ and all $X=(X_t)_{t\ge 0}\in D$
with $\sup_{t\le T}|X_t-\phi_t|\le\delta_{T,\phi}$ we have for
some $C_1>0$ and $C_2>0$
\[
C_1|\Lambda_{T,\phi}(y;t,X)|\ge g_{T,\phi}(y)\le
C_2(|\Lambda_{T,\phi}(y;t,\phi)| +1).
\]
Then condition {\bf IV} holds.
\end{lemma}

\medskip
{\it Proof.} First prove that
\begin{equation}
\int_0^T|\Lambda_{T,\phi}(\dot{\phi}_t;t,\phi)\dot{\phi}_t|dt<\infty.
\end{equation}
It suffices to verify (9.5) separately for
$\Lambda^+_{T,\phi}=\Lambda_{T,\phi} \vee 0$ and for
$\Lambda^-_{T,\phi}=-[\Lambda_{T,\phi}\wedge 0]$. By the
definition of $I_T(\phi)$
\begin{eqnarray}
&&
\int_0^T\big(\Lambda^+_{T,\phi}(\dot{\phi}_t;t,\phi)\dot{\phi}_t-
G(\Lambda^+_{T,\phi}(\dot{\phi}_t;t,\phi);t,\phi)\big)dt
\notag\\
&\le&\int_0^T\sup_\lambda(\lambda\dot{\phi}_t-G(\lambda;t,\phi))dt=I_T(\phi)
<\infty.
\end{eqnarray}
We saw in the proof of Lemma 9.1. (see (9.1) and (9.2)) that if
$\Lambda_{T,\phi}(\dot{\phi}_t;t,\phi)\ge 0$, then for all
$\lambda$ with $9\le\lambda\le
\Lambda^+_{T,\phi}(\dot{\phi}_t;t,\phi)$ for all $t\le T$
\[
\dot{\phi}_t\ge g(\lambda;t,\phi),
\]
and hence for any $m>1$
\begin{eqnarray}
\Lambda^+_{T,\phi}(\dot{\phi}_t;t,\phi)\dot{\phi}_t-
G(\Lambda^+_{T,\phi}(\dot{\phi}_t;t,\phi);t,\phi)&=&
\int_0^{\Lambda^+_{T,\phi}(\dot{\phi}_t;t,\phi)}(\dot{\phi}_t-g(\lambda;t,\phi))
d\lambda
\notag\\
&\ge&
\int_0^{\Lambda_{T,\phi}(\dot{\phi}_t;t,\phi)}(\dot{\phi}_t-g(\lambda;t,\phi))
d\lambda.
\end{eqnarray}
$g(\lambda;t,\phi)$ is increasing in $\lambda$ since (see (9.3))
$g'_\lambda(\lambda;t,\phi)\ge 0$ and so by (9.7)
\begin{eqnarray}
&& \Lambda^+_{T,\phi}(\dot{\phi}_t;t,\phi)\dot{\phi}_t-
G(\Lambda^+_{T,\phi}(\dot{\phi}_t;t,\phi);t,\phi)
\notag\\
&\ge&\frac{1}{m}\Lambda^+_{T,\phi}(\dot{\phi}_t;t,\phi)\Big(\dot{\phi}_t-
g\Big(\frac{1}{m}\Lambda^+_{T,\phi}(\dot{\phi}_t;t,\phi);t,\phi\Big)\Big).
\end{eqnarray}
Since by the definition of $\Lambda_{T,\phi}$, $\dot{\phi}_t=
g(\Lambda^+_{T,\phi}(\dot{\phi}_t;t,\phi);t,\phi)$, if
$\Lambda_{T,\phi}>0$, (9.2) implies that if
$\Lambda^+_{T,\phi}(\dot{\phi}_t;t, \phi)>0$ then
\begin{eqnarray}
&& \dot{\phi}_t-
g\Big(\frac{1}{m}\Lambda^+_{T,\phi}(\dot{\phi}_t;t,\phi);t,\phi\Big)
\notag\\
&=&\Big(1-\frac{1}{m}\Big)
\Lambda^+_{T,\phi}(\dot{\phi}_t;t,\phi)b^2(t,\phi)
\notag\\
&&+\int_Ef(t,\phi,u)\Big(\exp(\Lambda^+_{T,\phi}(\dot{\phi}_t;t,\phi)f((t,\phi,u)
)-1\Big)q(du)
\notag\\
&&-
\int_Ef(t,\phi,u)\Big(\exp\Big(\frac{1}{m}\Lambda^+_{T,\phi}(\dot{\phi}_t;t,\phi)
f((t,\phi,u)\Big)-1\Big)q(du).
\end{eqnarray}
Using the inequality $e^x-1\ge m(e^{x/m}-1), \ x\ge 0$, we have
\begin{eqnarray}
&&
\int_Ef^+(t,\phi,u)\Big(\exp(\Lambda^+_{T,\phi}(\dot{\phi}_t;t,\phi)f((t,\phi,u)
)-1\Big)q(du)
\notag\\
&&-
\int_Ef^+(t,\phi,u)\Big(\exp\Big(\frac{1}{m}\Lambda^+_{T,\phi}(\dot{\phi}_t;t,\phi)
f((t,\phi,u)\Big)-1\Big)q(du)
\notag\\
&\ge&(m-1)\int_Ef^+(t,\phi,u)\Big(\exp\Big(\frac{1}{m}
\Lambda^+_{T,\phi}(\dot{\phi}_t;t,\phi)f((t,\phi,u)\Big)-1\Big)q(du),
\end{eqnarray}
and since $e^{x/m}\ge e^x, \ x\le 0,$ so
\begin{eqnarray}
&&
\int_Ef^-(t,\phi,u)\Big(\exp(\Lambda^+_{T,\phi}(\dot{\phi}_t;t,\phi)f((t,\phi,u)
)-1\Big)q(du)
\notag\\
&&-\int_Ef^-(t,\phi,u)\Big(\exp\Big(\frac{1}{m}
\Lambda^+_{T,\phi}(\dot{\phi}_t;t,\phi)f((t,\phi,u)\Big)-1\Big)q(du)\le
0.
\end{eqnarray}
Now we substitute (9.10) and (9.11) in turn into (9.8) to obtain
\begin{eqnarray}
&& \Lambda^+_{T,\phi}(\dot{\phi}_t;t,\phi)\dot{\phi}_t-
G(\Lambda^+_{T,\phi}(\dot{\phi}_t;t,\phi);t,\phi)
\notag\\
&\ge&\frac{m-1}{m}\Lambda^+_{T,\phi}(\dot{\phi}_t;t,\phi)\Big[\frac{1}{m}
\Lambda^+_{T,\phi}(\dot{\phi}_t;t,\phi)b^2(t,\phi)
\notag\\
&&+\int_Ef^+(t,\phi,u)\Big(\exp\Big(\frac{1}{m}
\Lambda^+_{T,\phi}(\dot{\phi}_t;t,\phi)f((t,\phi,u)\Big)-1\Big)q(du).
\end{eqnarray}
(9.12) and (9.6) imply that
\begin{equation}
\int_0^Tb^2(t,\phi)(\Lambda^+_{T,\phi}(\dot{\phi}_t;t,\phi))^2dt<\infty
\end{equation}
\begin{eqnarray}
&& \int_0^T\Big[\int_Ef^+(t,\phi,u)\Big(\exp\Big(\frac{1}{m}
\Lambda^+_{T,\phi}(\dot{\phi}_t;t,\phi)f((t,\phi,u)\Big)-1\Big)q(du)\Big]
\notag\\
&&\times\Lambda^+_{T,\phi}(\dot{\phi}_t;t,\phi)dt<\infty.
\end{eqnarray}
(9.6), (9.8), and (9.13) and (9.14) lead in view of (9.2) lead to
\begin{eqnarray}
&&
\int_0^T\Lambda^+_{T,\phi}(\dot{\phi}_t;t,\phi)\Big[\dot{\phi}_t-a(t,\phi)+
\int_Ef^-(t,\phi,u)
\notag\\
&&\times\Big( \exp\Big(-\frac{1}{m}
\Lambda^+_{T,\phi}(\dot{\phi}_t;t,\phi)f^-((t,\phi,u)\Big)-1\Big)q(du)\Big]dt
<\infty.
\end{eqnarray}
By condition {\bf I} for $f$ and the inequality $1-e^{-x}\le x$,
we have
\begin{eqnarray}
&& \int_Ef^-(t,\phi,u)\Big[1-\exp\Big(-\frac{1}{m}
\Lambda^+_{T,\phi}(\dot{\phi}_t;t,\phi)f^-((t,\phi,u)\Big)\Big]q(du)
\notag\\
&\le&\frac{1}{m}\Lambda^+_{T,\phi}(\dot{\phi}_t;t,\phi)
\int_E(f^-((t,\phi,u))^2q(du)
\notag\\
&\le&\frac{1}{m}\Lambda^+_{T,\phi}(\dot{\phi}_t;t,\phi)
(1+\phi^*_T)^2l^2_T\int_Eh^2(u)q(du).
\end{eqnarray}
So far $m$ is arbitrary, we take it such that
\[
m\ge 2C_1C(1+\phi^*_T)^2l^2_T\int_Eh^2(u)q(du),
\]
with $C_1$ and $C$ from the assumptions of the lemma. Then (9.16)
and assumptions of the Lemma yield
\[
\int_Ef^-(t,\phi,u)\Big[1-\exp\Big(-\frac{1}{m}
\Lambda^+_{T,\phi}(\dot{\phi}_t;t,\phi)f^-((t,\phi,u)\Big)\Big]q(du)\le
1+\frac{|\dot{\phi}_t|}{2}.
\]
This together with (9.15) implies (note that the integrand in
(9.15) is nonnegative)
\begin{equation}
\int_0^T\Lambda^+_{T,\phi}(\dot{\phi}_t;t,\phi)I(\dot{\phi}_t>0)\Big(\dot{\phi}_t
-a(t,\phi)-\frac{1+\dot{\phi}_t}{2}\Big)dt<\infty.
\end{equation}
As by condition {\bf I}, $|a(t,\phi)|\le l_T(1+\phi^*_T), \ t\le
T$, so
\[
\Big(\dot{\phi}_t-a(t,\phi)-\frac{1+\dot{\phi}_t}{2}\Big)I(\dot{\phi}_t\ge
4l_T(1+\phi^*_T)+2)\ge \frac{\dot{\phi}_t}{4}I(\dot{\phi}_t\ge
4l_T(1+\phi^*_T)+2)
\]
and hence by (9.17)
\begin{equation}
\int_0^T\Lambda^+_{T,\phi}(\dot{\phi}_t;t,\phi)|\dot{\phi}_t|I(\dot{\phi}_t\ge
4l_T(1+\phi^*_T)+2)dt<\infty.
\end{equation}
Obviously $\Lambda^+_{T,\phi}(\dot{\phi}_t;t,\phi)=0$ if
$\dot{phi}_t\le -l_T(1+ \phi^*_T)$, also by the conditions of the
lemma $|\Lambda^+_{T,\phi}(\dot{\phi}_t;t,\phi)|\le
C_1C(1+|\dot{\phi}_t|)$, and so by (9.18)
\[
\int_0^T\Lambda^+_{T,\phi}(\dot{\phi}_t;t,\phi)|\dot{\phi}_t|dt<\infty.
\]
Argument for $\Lambda^-_{T,\phi}(\dot{\phi}_t;t,\phi)$ is similar.
(9.5) is proved.

We end the proof. Let $\phi^N=(\phi^N_t)_{t\ge 0}$ be defined in
condition {\bf IV}. Since $\phi\in D$ is absolutely continuous so
$\int_0^T|\dot{\phi}_t|dt <\infty$ and we may regard that
\begin{equation}
\sup_{t\le T}|\phi^N_t-\phi_t|\le
\int_0^T|\dot{\phi}_t|I(|\dot{\phi}_t|>N)dt\le \delta_\phi.
\end{equation}
By Remark to condition {\bf III} (see (2.1))
\begin{eqnarray}
H(\dot{\phi}_t;t,\phi^N)&=&\Lambda_{T,\phi}(\dot{\phi}_t;t,\phi^N)\dot{\phi}_t-
G(\Lambda_{T,\phi}(\dot{\phi}_t;t,\phi^N);t,\phi^N)
\notag\\
&\le&|\Lambda_{T,\phi}(\dot{\phi}_t;t,\phi^N)\dot{\phi}_t|+
|\Lambda_{T,\phi}(\dot{\phi}_t;t,\phi^N)a(t,\phi^N)|.
\end{eqnarray}
In view of conditions {\bf I} and {\bf II}, and (9.19) we may well
assume that for large $N$
\[
\sup_{t\le T}a(t,\phi^N)|\le l_T(2+\phi^*_T),
\]
and then (9.20) yields for almost all $t\le T$ that
\[
H(\dot{\phi}_t;t,\phi^N)\le
C'|\Lambda_{T,\phi}(\dot{\phi}_t;t,\phi^N)| (1+|\dot{\phi}_t|)
\]
for some $C'>0$. As the assumptions of the lemma and (9.19) imply
that
$$
|\Lambda_{T,\phi}(\dot{\phi}_t;t,\phi^N)|\le
C^{-1}_1g_{T,\phi}(\dot{\phi}_t),
$$
we thus obtain
\begin{equation}
H(\dot{\phi}_t;t,\phi^N)\le CC^{-1}_1g_{T,\phi}(\dot{\phi}_t)
(1+|\dot{\phi}_t|).
\end{equation}
By the assumption of the lemma and (9.5)
\[
\int_0^Tg_{T,\phi}(\dot{\phi}_t)(1+|\dot{\phi}_t|)dt<\infty,
\]
and so (9.21) gives condition {\bf IV}.

\medskip
{\it Proof of Theorem 9.1.} Condition {\bf III} holds by Lemma
9.1. We now prove that conditions of Lemma 9.2 hold.

Assume first that $f^+(t,X,u)$ and $f^-(t,X,u)$ are nondegenerate
with respect to $q(du)$ uniformly on $[0,T]$ in a neighborhood of
$\phi$.

If $y\ge l_T(2+\phi^*_T)$, then by conditions {\bf I} and {\bf II}
there exists $\delta\le \delta_{T,\phi}$ such that $sup_{t\le
T}|a(t,X)|\le y$ for all $X= (X_t)_{t\ge 0}\in D$ with $sup_{t\le
T}|X_t-\phi_t|\le \delta$. The by definition of $\Lambda_{T,\phi}$
for these $X$, $\Lambda_{T,\phi}(y;t,X)\ge 0$ and hence for any
$\varepsilon>0$
\[
y\ge
-l_T(2+\phi^*_T)+[\exp(\Lambda_{T,\phi}(y;t,X)\varepsilon)-1]\int_E
f^+(t,X,u)I(f^+(t,X,u)\ge \varepsilon)q(du),
\]
which implies that there exists $r>0$ such that
\begin{equation}
\Lambda_{T,\phi}(y;t,X)\le r\ln(1+y)
\end{equation}
for all $t\le T$, all $X=(X_t)_{t\ge 0}\in D$ with $\sup_t{\le
T}|X_t-\phi_t|\le \delta$ and $y\ge l_T(2+\phi^*_T)$.

The nondegeneracy of $f^-(t,X,u)$ provides us with an estimate
similar to (9.22) for $y\le -l_T(2+\phi^*_T)$.

Hence we can find $r_1$ such that for all $y\in R$, all
$X=(X_t)_{t\ge 0}$ with
\[
\sup_{t\le T}|X_t-\phi_t|\le\delta
\]
and all $t\le T$ we have
\begin{equation}
|\Lambda_{T,\phi}(y;t,X)|\le r_1\ln(1+|y|).
\end{equation}

On the other hand, the definition of $\Lambda_{T,\phi}$ and
condition {\bf I} imply that for all $t\le T$
\begin{eqnarray}
|y|&\le&l_T(1+\phi^*_T)+|\Lambda_{T,\phi}(y;t,X)|l^2_T(1+\phi^*_T)^2
\notag\\
&&+(|\Lambda_{T,\phi}(y;t,X)|+1)l^2_T(1+\phi^*_T)^2\int_Eh^2(u)q(du)
\notag\\
&&+K((|\Lambda_{T,\phi}(y;t,X)|+1)l_T(1+\phi^*_T)
\end{eqnarray}
($K(\lambda)$ is defined in condition {\bf I}). By assumptions of
the theorem there exist $A_1$ and $A_2$ such that $K(\lambda)\le
A_1\exp(\lambda A_2)$ for all $\lambda\ge 0$, and it then follows
in view of (9.24) that for some $r_1>0$ which depends on $T$ and
$\phi$ we have
\begin{equation}
|\Lambda_{T,\phi}(y;t,X)|\ge r_1\ln(1+|y|)
\end{equation}
for all $t\le T$.

(9.23) and (9.25) show that the conditions of Lemma 9.2 are met by
$g_{T,\phi}(y)=\ln(1+|y|)$.

Now let $f^+(t,X,u)$ and $f^-(t,X,u)$ be degenerate with respect
to $q(du)$ uniformly on $[0,T]$ in a neighborhood of $\phi$. Then
by (9.1) and (9.2) and the definition of $\Lambda_{T,\phi}$ we
have for all $t\le T$ and all $X=(X_t)_{t\ge 0}\in D$ with
$\sup_{t\le T}|X_t-\phi_t|\le \delta$
\[
y\ge-l_T(2+\phi^*_T)+\Lambda_{T,\phi}(y;t,X)b^2(t,X)
\]
if $y\ge l_T(+\phi^*_T)$, and
\[
y\le l_T(2+\phi^*_T)-\Lambda_{T,\phi}(y;t,X)b^2(t,X)
\]

if $y\le -l_T(+\phi^*_T)$.

So if $b(t,X)$ is uniformly nondegenerate, than
\[
|\Lambda_{T,\phi}(y;t,X)|\le L(1+|y|)
\]
for all the above $t$ and $X$ and all $y\in R$ and some $L>0$. As,
obviously
\[
|y|\le l_T(1+\phi^*_T)+|\Lambda_{T,\phi}(y;t,X)|b^2(t,X),
\]
so conditions of Lemma 9.3 hold with $g_{T,\phi}(y)=1+|y|$.

Finally,, it is easily follows from the above that if $f^+(t,X,u)$
is degenerate with respect  to $q(du)$ uniformly on $[0,T]$ in a
neighborhood of $\phi$, $f^-(t,X,u)$ and $b(t,X)$ are both
nondegenerate then we can take
\[
g_{T,\phi}(y)=
  \begin{cases}
    1+y & y\ge 0, \\
    1+\ln(1-y), & y\le 0,
  \end{cases}
\]
and if degenerate is $f^-(t,X,u)$ and $f^+(t,X,u)$ and $b(t,X)$
are nondegenerate then conditions of Lemma 9.2 are met by
\[
g_{T,\phi}(y)=
  \begin{cases}
    1+\ln(1+y), & y\ge 0, \\
    1-y, & y\le 0.
  \end{cases}
\]

\end{document}